\documentclass[11pt]{amsart}
\usepackage{mathrsfs}  
\usepackage[backref,colorlinks]{hyperref} 
\usepackage{amsmath,amssymb,amsthm,bm}
\usepackage[foot]{amsaddr}
\usepackage[abbrev,msc-links]{amsrefs}   
\usepackage[utf8]{inputenc}
\usepackage{tikz}
\usepackage{enumerate}
\usepackage{enumitem}
\usepackage{graphicx}
\usepackage{euflag}

\usepackage[normalem]{ulem}

\usepackage{mathtools}
\DeclarePairedDelimiter\ceil{\lceil}{\rceil}
\DeclarePairedDelimiter\floor{\lfloor}{\rfloor}

\usepackage{lineno}
\theoremstyle{plain}
\newtheorem{thm}{Theorem}[section]

\newtheorem{clm}[thm]{Claim}
\newtheorem{cor}[thm]{Corollary}
\newtheorem{lem}[thm]{Lemma}

\newtheorem*{clm*}{Claim}

\theoremstyle{definition}
\newtheorem{rem}[thm]{Remark}
\newtheorem{dfn}[thm]{Definition}
\newtheorem{exmp}[thm]{Example}

\numberwithin{equation}{section}

\normalsize                              
\newenvironment{claimproof}[1]{\par\noindent{\it Proof of Claim:}\space#1}{\hfill $\blacksquare$}

\def\COMMENT#1{}
\let\COMMENT=\footnote

\let\polishlcross=\l
\def\l{\ifmmode\ell\else\polishlcross\fi}

\newcommand{\bE}{\mathbb{E}}

\newcommand{\bP}{\mathbb{P}}
\def\NN{\mathbb N}

\def\cA{{\mathcal A}}
\def\cH{{\mathcal H}}

\def\cF{{\mathcal F}}
\def\cB{{\mathcal B}}
\def\cC{{\mathcal C}}

\def\cP{{\mathcal P}}
\def\cY{{\mathcal Y}}
\def\cR{{\mathcal R}}
\def\cQ{{\mathcal Q}}
\def\cE{{\mathcal E}}
\def\cX{{\mathcal X}}

\def\cI{{\mathcal I}}

\def\eps{\varepsilon}
\def\phi{\varphi}

\def\vecsign{\mathchar"017E}
\def\dvecsign{\smash{\stackon[-1.95pt]{\vecsign}{\rotatebox{180}{$\vecsign$}}}}
\def\dvec#1{\def\useanchorwidth{T}\stackon[-4.2pt]{#1}{\,\dvecsign}}
\usepackage{stackengine}
\stackMath

\DeclareMathAlphabet\bfc{OMS}{cmsy}{b}{n}
\DeclareMathAlphabet{\pzc}{OT1}{pzc}{m}{it}

\DeclareMathAlphabet\bc{OMS}{cmsy}{b}{n}

\title[A rainbow Dirac theorem for loose Hamilton cycles in hypergraphs]{A rainbow Dirac theorem for loose Hamilton cycles in hypergraphs}
 \author{Amarja Kathapurkar $^{1,\ast}$}
 \author{Patrick Morris $^{2,\dagger}$}
 \author{Guillem Perarnau $^{2,3,\ddagger}$}

 \address{$^1$ University of Birmingham, United Kingdom.}
 \address{$^2$ Universitat Polit\`ecnica de Catalunya (UPC), Barcelona, Spain.}
\address{$^3$ Centre de Recerca Matem\`atica, Bellaterra, Spain.}

\thanks{$^\ast$ Research supported by EPSRC Research grant EP/R034389/1. }

\thanks{$^\dagger$ Research supported  by the Deutsche Forschungsgemeinschaft (DFG, German
Research Foundation) Walter Benjamin programme - project number 504502205 and by the European Union's Horizon Europe   Marie Sk{\l}odowska-Curie grant RAND-COMB-DESIGN - project number
101106032 {\euflag}.}

\thanks{$^\ddagger$ Research
supported by the grants RED2022-134947-T, PID2023-147202NB-I00, PCI2024-155080-2 and the Programme Severo Ochoa y Mar\'ia de Maeztu por Centros y Unidades de Excelencia en I\&D (CEX2020-001084-M), all of them funded by MICIU/AEI/10.13039/501100011033.}

\email{amarja.kathapurkar@gmail.com, pmorrismaths@gmail.com,  guillem.perarnau@upc.edu}

\date{\today}

\begin{document}

\begin{abstract}
A meta-conjecture of Coulson, Keevash, Perarnau and Yepremyan \cite{ckpy2020rainbow} states that above the extremal threshold for a given spanning structure in a (hyper-)graph, one can find a rainbow version of that spanning structure in any suitably bounded colouring of the host (hyper-)graph. We solve one of the most pertinent outstanding cases of this conjecture,  by showing that for any $1\leq j\leq k-1$, if $G$ is a $k$-uniform hypergraph above the $j$-degree threshold for a loose Hamilton cycle, then any globally bounded colouring of $G$ contains a rainbow loose Hamilton cycle.
\end{abstract}
\maketitle 
\section{Introduction}
A seminal theorem of Dirac \cite{dirac1952some} states that any $n$-vertex graph $G$ with $\delta(G)\geq n/2$ contains a Hamilton cycle. This inspired many further results exploring the optimal minimum degree conditions for certain spanning structures in a host (hyper-)graph and this area, sometimes referred to as `Dirac theory', is a cornerstone of modern extremal combinatorics. The field has flourished in recent decades due to powerful tools being developed to tackle these questions, such as the regularity method  and absorption. In graphs, this has led to a deep understanding of the full picture. 

In hypergraphs, the situation is considerably more complex. This is, in part, due to the various ways in which one can generalise the graph case. For example, when generalising Dirac's theorem to hypergraphs, one has a range of choices as to which minimum degree condition is considered and what type of Hamilton cycle is desired. Indeed, for a $k$-uniform hypergraph $G$ ($k$-graph for short), one can consider 
\[\delta_j(G):=\min\left\{|\{e\in E(G):T\subset e\}|:T\in \binom{V(G)}{j}\right\},\]
 for $1\leq j \leq k-1$.  Likewise with Hamilton cycles, one can consider a cyclic ordering of the vertices of $G$ and require that each edge of the Hamilton cycle occupies $k$ consecutive vertices in the ordering and every pair of consecutive edges intersect in precisely $\ell$ vertices for some $1\leq \ell\leq k-1$. Such a Hamilton cycle is called a \emph{Hamilton $\ell$-cycle} and when $\ell=1$, we refer to it as a \emph{loose} cycle, whilst the case $\ell=k-1$ is referred to as a \emph{tight} cycle. Note that if an $n$-vertex $k$-graph $G$ has a Hamilton $\ell$-cycle, then one necessarily has that $(k-\ell)|n$. 

In hypergraphs our understanding of minimum degree thresholds is far from complete, even in simple cases such as perfect matchings (see e.g. \cite{zhao2016recent}).  
 For Hamilton cycles, we make the following definition to capture the relevant extremal minimum degree conditions. 

\begin{dfn}[Hypergraph Hamilton cycle thresholds]\label{def:threshold}
For $1\leq j, \ell\leq k-1$, we denote by $\delta_j^k(\ell)$ the smallest real number 
$\delta\geq 0$ such that for  any $\eps>0$ and  $n\in (k-\ell)\NN$ sufficiently large, any $n$-vertex $k$-graph $G$ with $\delta_j(G)\geq (\delta+\eps)n^{k-j}$ contains a Hamilton $\ell$-cycle.     
\end{dfn}

Determining these threshold values has proven to be a considerable challenge. For example, in the case of loose Hamilton cycles (when $\ell=1$),  a case  that has been of particular interest, we only know the location of $\delta_j^k(1)$ when $j=k-2$ or $j=k-1$. For $j=k-1$, the so-called \emph{codegree threshold}, this was established by K\"uhn and Osthus \cite{kuhn2006loose} when $k=3$ and in general by Keevash, K\"uhn, Mycroft and Osthus \cite{keevash2011loose} and, independently, H\`an and Schacht \cite{han2010dirac}. For $j=k-2$, this was again solved first in the $3$-uniform case, by Bu\ss, H\`an and Schacht \cite{buss2013minimum}, and then for general $k$ by Bastos, Mota, Schacht, Schnitzer and Schulenburg \cite{de2017loose}. 
For other values of $\ell$, the results are also sporadic, we refer the reader to the survey~\cite{zhao2016recent} for a comprehensive review.



Beyond determining extremal degree thresholds, another fruitful line of research in Dirac theory has been to show \emph{robustness} for these thresholds. 
 The aim here is to show that when a (hyper-)graph is above the minimum degree threshold  (we will informally refer to such (hyper-)graphs as being `Dirac') with respect to a given spanning structure, the Dirac (hyper-)graph is in fact \emph{robust} with respect to containing that spanning structure. This includes considerations such as the existence of many spanning structures and resilience to random edge deletions.  We refer to the nice survey of Sudakov \cite{sudakov2017robustness} and the references therein for an introduction to these results in the context of Dirac's condition for Hamilton cycles. 
In this paper, we  consider a notion of robustness  related to  finding  \emph{rainbow} spanning structures in any bounded edge colouring of the Dirac (hyper-)graph. This is motivated by the classical study of rainbow spanning structures in certain colourings of graphs. 

\subsection{Rainbow spanning structures} 
A subgraph $H$ of an edge coloured graph $G$ is  \emph{rainbow} if each of the edges of $H$ is a different colour. Rainbow subgraphs appeared early on in combinatorics via connections to design theory with several beautiful conjectures being posed. 
Perhaps the most famous such conjecture, known as the Ryser-Brualdi-Stein conjecture \cite{brualdi1991combinatorial,ryser1967neuere,stein1975transversals} states that every $n\times n$ Latin square has a transversal of size at least $n-1$  and one of size $n$ when $n$ is odd. The first part of this (establishing the existence of transversals of size $n-1$) has only recently been solved  by Montgomery \cite{montgomery2023proof}. Translating to colourings of graphs, the Ryser-Brualdi-Stein conjecture is equivalent to the assertion that one can always find a rainbow matching of size $n-1$ (or $n$ when $n$ is odd) in any properly coloured $K_{n,n}$ with $n$ colours. 

From a graph theoretic perspective one can ask more generally what conditions on a colouring of a host graph guarantee the existence of a rainbow (almost) spanning structure of interest. 
Indeed 
a colouring being proper is equivalent to saying that the colouring is \emph{$t$-locally bounded}  (at each vertex we see every colour at most $t$ times) with $t=1$.  
One can also then consider \emph{globally bounded} conditions where we instead bound the size of each colour class in the whole (hyper-)graph.  The first result of this kind was due to Erd\H{o}s and Spencer \cite{erdos1991lopsided}, who showed that any colouring of $K_{n,n}$ with at most $\tfrac{n}{16}$ edges of each colour contains a rainbow perfect matching.
Another early example of interest was due to  Erd\H{o}s and Stein (see \cite{erdos1983some}) who asked whether there is some constant $c>0$ such that any  colouring of $K_n$  with at most $cn$ edges of each colour contains a rainbow Hamilton cycle. This was then explicitly conjectured by Hahn and Thomassen \cite{hahn1986path} and, after several results towards the conjecture,  was solved by Albert, Frieze and Reed \cite{albert1995multicoloured}. A generalisation to hypergraph Hamilton cycles was then given by Dudek, Frieze and Ruci\'nski \cite{dudek2012rainbow}. There has been a wealth of similar results studying different spanning structures. 


\subsection{Rainbow structures in Dirac (hyper-)graphs}
The  majority of results concerning rainbow spanning substructures in  bounded (and proper) colourings have focused on the case where the host graph is a complete (hyper-)graph or complete bipartite graph. When considering other possible host graphs,  Dirac graphs arise naturally.  This perspective was first considered by Cano, Perarnau and Serra \cite{cano2017rainbow} who showed that one can find a rainbow Hamilton cycle in any globally $o(n)$-bounded colouring of $G$ when $G$ is either an $n$-vertex graph or a balanced bipartite graph with $n$ vertices in each part, and such that $G$ has minimum degree $\delta(G)\geq (1+o(1))\tfrac{n}{2}$. The asymptotic minimum degree condition was then replaced to give an exact minimum degree condition    $\delta(G)\geq \tfrac{n}{2}$ by Coulson and Perarnau, first in the bipartite case \cite{coulson2019rainbow} and then in the non-bipartite case \cite{coulson2020rainbow} as in Dirac's original theorem. These results thus give evidence of robustness for the extremal thresholds for Hamilton cycles. Note also that in the bipartite case, these results can be seen as a direct strengthening of the result of Erd\H{o}s and Spencer \cite{erdos1991lopsided}, allowing for host graphs that are not complete (at the expense of a  worse constant for the boundedness). 

Further examples of these types of results came from Coulson, Keevash, Perarnau and Yepremyan \cite{ckpy2020rainbow} who proved that (asymptotically) above the minimum degree for a given (hyper-)graph $F$-factor, one finds a rainbow $F$-factor in any suitably bounded colouring, and from Glock and Joos \cite{glock2020rainbow} who gave a rainbow version of the famous blow-up lemma \cite{komlos1997blow}, which allowed them to give results of this flavour in considerable generality for graphs, in particular providing a rainbow version of the bandwidth theorem \cite{bottcher2009proof}.  We  remark that a nice feature of the work of \cite{ckpy2020rainbow} is that they could establish such a result, even in cases where the minimum degree threshold has not yet been determined. 

All of these results provide evidence of a general phenomenon and caused Coulson, Keevash, Perarnau and Yepremyan \cite{ckpy2020rainbow} to explicitly give the ``meta-conjecture" that once one is above the extremal threshold for a given spanning structure, rainbow copies of that structure can be found in any suitably bounded colouring of the Dirac graph. Our main result provides further evidence for this conjecture by establishing that this is the case for loose Hamilton cycles in hypergraphs.

\begin{thm}[Main Theorem] \label{thm:main} Let $1\leq j< k\in \NN $ and let $\eps>0$ be arbitrary. Then there exists $\mu>0$ such that for any sufficiently large $n\in (k-1)\NN$, we have that if $G$ is an $n$-vertex $k$-graph with $\delta_j(G)\geq (\delta_j^k(1)+\eps)n^{k-j}$, the following holds. For any colouring $\chi:E(G)\rightarrow \NN$  of $G$ such that \begin{enumerate}[label={\textbf{(B)}}]
    \item \label{item:Global}  $|\chi^{-1}(c)|\leq \mu n^{k-1}$ for all $c\in \NN$,
\end{enumerate} 
we have that there is a rainbow loose Hamilton 
cycle in $G$ coloured by $\chi$. 
    \end{thm}

Theorem \ref{thm:main}  provides a first generalisation of the result of Coulson and Perarnau \cite{coulson2020rainbow} to the hypergraph setting and strengthens the previously mentioned work of Dudek, Frieze and Ruci\'nski  \cite{dudek2012rainbow} who proved the statement in the case that the host hypergraph $G$ is complete.
    Note that the bound in~\ref{item:Global} is tight up to the choice of $\mu$ as some bound of the order of $n^{k-1}$ is required in order to have sufficiently many colours available for a rainbow loose Hamilton cycle.  Note also that, as with the previous results on $F$-factors in hypergraphs  due to Coulson, Keevash, Perarnau and Yepremyan \cite{ckpy2020rainbow}, our result works even for the values of $j$ where the threshold value $\delta_j^k(1)$ is not yet known.

    We remark that Theorem \ref{thm:main} and the aforementioned result \cite{ckpy2020rainbow} on $F$-factors are the only results of this kind, proving robustness via rainbow copies in bounded colourings, in the realm of hypergraphs. It is not surprising that progress has been much slower when dealing with higher uniformities as many of the tools used in the graph setting \cite{coulson2019rainbow,glock2020rainbow}, in particular the regularity method, face barriers for hypergraphs in Dirac theory, as discussed above. Our proof of Theorem \ref{thm:main} follows the same scheme as that of \cite{ckpy2020rainbow} for $F$-factors, appealing to the local lemma, and using random samples to perform switches. The main barrier in applying these methods in the setting of Hamilton cycles is that the structure we are trying to find is \emph{connected}. This prevents us from being able to locally adapt a copy of our desired spanning structure independently from the structure as a whole. These local adaptations are at the core of the switching method used in \cite{ckpy2020rainbow} and, at a high level, we overcome this obstacle by using absorption techniques to piece back together loose Hamilton cycles after local alterations. Making this work requires new ideas to overcome technical hurdles at many parts of the proof and we also require recent work of Alvarado, Kohayakawa, Lang, Mota and Stagni \cite{alvarado2023resilience} who showed the existence of an absorption scheme when the location of the threshold is not known. Another technical contribution of our work is to remove the need for a local bound on the colouring. Indeed, the result of Coulson, Keevash, Perarnau and Yepremyan \cite{ckpy2020rainbow} found $F$-factors under the assumption \ref{item:Global} as well as an additional constraint that every $(k-1)$-set lies in at most $\mu n$ edges of each colour. This local boundedness condition is superfluous when $k=2$ (the graph case) but it is not implied by the global condition \ref{item:Global} when $k\geq 3$. The authors of \cite{ckpy2020rainbow} justify this extra condition by showing that for certain $F$, for example when $F$ is a clique with more than one edge, the condition is in fact \emph{necessary}. However in the setting of loose Hamilton cycles, Theorem \ref{thm:main} shows that the global condition \ref{item:Global} suffices.  Within the proof scheme of \cite{ckpy2020rainbow} and our work here, removing the need for the local condition requires new ideas to maintain the existence of a rainbow structure. We develop  a novel strategy using random partitions and deterministic adjustments to the partition and we remark that these ideas  can be used to strengthen the result of \cite{ckpy2020rainbow} to be able to drop the local boundedness condition on the colouring to get $F$-factors when $F$ is \emph{linear}, that is, when any pair of edges in $F$ intersect in at most one vertex. We discuss this further in the concluding remarks, Section \ref{sec:conclude}.

    \subsection{Organisation} In the next section, we collect some notation and terminology as well as some probabilistic tools that we will use throughout. In Section \ref{sec:overview}, we then provide a proof overview and reduce the proof of Theorem \ref{thm:main} to three main lemmas which are proven in Sections \ref{sec:lll}, \ref{sec:absorption} and \ref{sec:random sample}. Finally in Section \ref{sec:conclude}, we discuss potential future directions of research.

\

\section{Preliminaries} \label{sec:prelims}

\subsection{Notation and terminology} \label{sec:notation}

For $k$-graphs $H$ and $G$, we use the notation $H\subset G$ to indicate that $H$ is a \emph{subhypergraph} (or \emph{subgraph} for short) of $G$. 
A \emph{loose path} (or simply path, for short) of \emph{length} $t$ in a $k$-graph $G$ is a subgraph $P$ of $G$ defined by a collection of $t$ edges $e_1,\ldots,e_t\in E(G)$ such that there is some ordering $v_1,\ldots, v_{t(k-1)+1}$ of $V(P)\subset V(G)$ with $e_i=\{v_{1+(i-1)(k-1)}, \ldots, v_{k+(i-1)(k-1)}\}$ for $i=1,\ldots,t$. The vertices $V^E(P):=\{v_1,v_{t(k-1)+1}\}$ are referred to as the \emph{endvertices} of the path. 

A loose Hamilton cycle $H$ comes with a cyclic order on $E(H)$ such that the only intersecting edges of $H$ are pairs of consecutive edges in the ordering (which intersect in exactly one vertex). For each loose Hamilton cycle $H$ of $G$, we will act under the convention that there is a fixed orientation of this order of $E(H)$ and will refer to an \emph{increasing} path, of length $t$ say, from an edge $e\in E(H)$ to be the path defined by taking $e$ and the $t-1$ edges that follow $e$ when traversing the loose Hamilton cycle according to the orientated order on $E(H)$. 

For a $k$-graph $G$ and $S\subset V(G)$ with $|S|\leq k$, we define the  set $E_G(S)\subseteq E(G)$ of edges \emph{incident to $S$} to be 
\[E_G(S):=\{e\in E(G):S\subseteq e\},\]
and we define the 
\emph{degree} of $S$ in $G$ to be 
\[\deg_G(S)=|E_G(S)|.\]
We write simply $E(S)$ and  $\deg(S)$ if the $k$-graph $G$ is clear from context. 
We will also need neighbourhoods and degrees defined relative to vertex subsets. 
For a  $k$-graph $G$ and  vertex subsets  $S, W\subset V(G)$ with $|S|\leq k$, we define the  set $E_{G}(S;W)\subseteq E(G)$  to be 
\[E_{G}(S;W):=\{e\in E(G):S\subseteq e, e\setminus S\subset W \},\]
and we define the 
\emph{degree} of $S$ in $G$ relative to $W$ to be 
\[\deg_{G}(S;W)=|E_{G}(S;W)|.\]
We write simply $E(S;W)$ and  $\deg(S;W)$ if the $k$-graph $G$ is clear from context.

\subsection{Probabilistic tools} \label{sec:tools}

In this section we collect various probabilistic tools that
will be used in the proof of Theorem~\ref{thm:main}.
We start with a general version of the (lopsided) local lemma, see \cite[Lemma 5.1.1]{alon16} and the remarks following the lemma.

\begin{lem}[The lopsided local lemma] \label{lem:llll}
Suppose $\cA$ is a set of events in a finite probability space, $\Gamma$ is a graph with $V(\Gamma)=\cA$ and there are real numbers $\{x_A:A\in \cA\}$ with $0\leq x_A <1$ for all $A\in \cA$. Suppose further that for any $A\in \cA$ and subset  $\mathcal{A}' \subseteq \mathcal{A}\setminus \{A\}$
such that $AA' \notin E(\Gamma)$ for all $A' \in \mathcal{A}'$ and $\bP[\cap_{A' \in \mathcal{A}'} \overline{A'}]>0$, we have that 
 \begin{linenomath}
\begin{equation} \label{eq:lll_cond_prob}
\bP[A\vert \cap_{A' \in \mathcal{A}'} 
\overline{A'}]\leq x_A
\prod_{B:AB\in E(\Gamma)} (1-x_B). 
\end{equation}
\end{linenomath}
Then with positive probability none of the events $A\in \cA$ occur. 
\end{lem}

\begin{rem} The fact that we only have a condition for families $\cA'$ for which  $\bP[\cap_{A' \in \mathcal{A}'} \overline{A'}]>0$ is not mentioned in \cite[Section 5.1]{alon16} (or elsewhere in the literature). However such a condition is necessary for the probability  in \eqref{eq:lll_cond_prob} to be well-defined. This omission presumably stems from the fact that Lemma~\ref{lem:llll} is in fact proven in \cite{alon16} under the stronger assumption that $A$ is mutually independent from all of the events $\{A':AA'\in E(\Gamma)\}$ and $\bP[A]$ satisfies the inequality in \eqref{eq:lll_cond_prob}. Therefore one can bypass the issue of having $\bP[\cap_{A' \in \mathcal{A}'} \overline{A'}]=0$ by using the mutual independence to deduce that $\bP[A\vert \cap_{A' \in \mathcal{A}'} 
\overline{A'}]=\bP[A]$ and bounding $\bP[A]$. Strictly speaking, there is still a potential issue here with  $\bP[A\vert \cap_{A' \in \mathcal{A}'} 
\overline{A'}]$ being well-defined but an inspection of the proof of \cite[Lemma 5.1.1]{alon16} shows that one can add an extra induction over the size of $|\cA'|$ which will show that $\bP[\cap_{A' \in \mathcal{A}'} \overline{A'}]>0$ for all $\cA'\subseteq \cA$ (and ultimately $\cA$ itself) and so imposing \eqref{eq:lll_cond_prob} only for $\cA'$ with $\bP[\cap_{A' \in \mathcal{A}'} \overline{A'}]>0$ is justified. 
\end{rem}

In fact, it will be easier for us to work with a slight modification of Lemma~\ref{lem:llll}, for which we make the following definition.

\begin{dfn} \label{def:dep}
Let $\cA$ be a set of events in a finite probability space.
Suppose $\Gamma$ is a graph with $V(\Gamma) = \mathcal{A}$ and $\textbf{p} \in [0,1/2)^{\mathcal{A}}$.
We call $\Gamma$ a \emph{$\textbf{p}$-dependency graph for $\cA$}
if for every $A \in \mathcal{A}$ and $\mathcal{A}' \subseteq \mathcal{A}\setminus \{A\}$
such that $AA' \notin E(\Gamma)$ for all $A' \in \mathcal{A}'$ and $\bP[\cap_{A' \in \mathcal{A}'} \overline{A'}]>0$,
we have $\bP[A \vert \cap_{A' \in \mathcal{A}'} \overline{A'}] \leq p_A$.
\end{dfn}

We now deduce the following corollary of Lemma~\ref{lem:llll} which we will use in applications. This result is already known in the literature but we add its proof for completeness.

\begin{cor} \label{L4}
Under the setting of Definition~\ref{def:dep},
if $\sum \{ p_{B} : AB \in E(\Gamma) \} \le 1/4$
for all $A \in \cA$ then with positive probability
none of the events in $\cA$ occur.
\end{cor}
\begin{proof}
For each $A\in \cA$, we set $x_A=2p_A$ and need to show that the condition \eqref{eq:lll_cond_prob} holds for such values. So fix some $A\in \cA$ and some 
$\mathcal{A}' \subseteq \mathcal{A}$
with $\bP[\cap_{A' \in \mathcal{A}'} \overline{A'}]>0$ and such that $AA' \notin E(\Gamma)$ for all $A' \in \mathcal{A}'$. Then 
\[x_A
\prod_{B:AB\in E(\Gamma)} (1-x_B)\geq x_A\left(1-\sum_{B:AB\in E(\Gamma)} x_B\right)=2p_A\left(1-2\sum_{B:AB\in E(\Gamma)} p_B\right)\geq p_A\geq \bP[A \vert \cap_{A' \in \mathcal{A}'} \overline{A'}],\]
as required.
\end{proof}

We will  frequently use  concentration inequalities for random variables. The first such inequality, Chernoff's inequality~\cite{Chernoff1952} (see also~\cite[Corollaries 2.3 and 2.4]{Janson2011}), deals with the case of binomial random variables.

\begin{thm}[Chernoff inequality] \label{thm:chernoff}
Let~$X$ be the sum of a set of mutually  independent Bernoulli random variables and let~$\lambda=\bE[X]$. Then for any~$0<\delta<\tfrac{3}{2}$, we have that 
\[\bP[X\geq (1+\delta)\lambda]\leq  e^{-\delta^2\lambda/3 } \hspace{2mm} \mbox{ and } \hspace{2mm} \bP[X\leq (1-\delta)\lambda] \leq  e^{-\delta^2\lambda/2 }.\]
 Moreover, if~$x\geq 7 \lambda$, then~$\bP[X\geq x]\leq e^{-x}$.
\end{thm}

The following well known bound gives a lower bound on the probability that a binomial random variable hits its expected value. The form we give here can be found in \cite[Lemma 3.6]{ckpy2020rainbow}.

\begin{lem} \label{lem:hit expected}
Let~$X$ be binomially distributed with parameters $({n},{p})$ and suppose that~$\lambda:=\bE[X]=np\in \NN$ is an integer with $\lambda=o(\sqrt{n})$. Then $\bP[X=\lambda]\geq 1/(4\sqrt{\lambda})$. 
\end{lem}

We will also be interested in concentration inequalities when dealing with random variables that are not mutually independent. Let  $\{X_i\}_{i\in \cI}$ be a family of random variables on a common probability space. A \emph{(strong) dependency graph} for $\{X_i\}_{i\in \cI}$ is any (2-uniform) graph $\Gamma$ with $V(\Gamma)=\cI$ and such that if $A$ and $B$ are disjoint subsets of $\cI$ and $\{e\in E(\Gamma):e\cap A, e\cap B\neq\emptyset\}=\emptyset$, then the families $\{X_i\}_{i\in A}$ and $\{X_i\}_{i\in B}$ are independent of each other.
The following inequality is a version of Suen's inequality \cite{suen1990correlation} given by Janson \cite{janson1998new}
 (see also~\cite[Theorem 2.23]{Janson2011}) and it  provides a bound for the lower tail in the case that there is a strong dependency graph which is sufficiently sparse.

\begin{lem}[Suen's inequality] \label{lem:suen}
Let $\{I_i\}_{i\in\cI}$ be a finite family of indicator random variables with some strong dependency graph $\Gamma$. Let~$X:=\sum_{i\in \cI} I_i$, ~$\lambda := \bE(X)$ and write $i\sim j$ if $ij\in E(\Gamma)$. Moreover,
let~$\Delta_{X} := \sum_{(i,j)\in \cI^2:i\sim j}\bE[I_i I_j]$ and $\delta_X=\max_{i\in \cI}\sum_{k\in \cI:k\sim i}\bP[I_k=1]$. 
Then  for any~$0\le t\le \lambda$,
 \begin{linenomath}
\begin{equation}
\bP[X\leq \lambda -t]\leq \exp \left ( -\min\left(\frac{t^2}{4(\lambda+\Delta_{X})}, \frac{t}{6\delta_X}\right)\right ). \label{eq:2}
\end{equation}
\end{linenomath}
\end{lem}

Finally we need the following inequality of Talagrand \cite{talagrand1995concentration}, which we give in the following form, see \cite[Section 10.1]{molloy2002graph}. 

\begin{lem}[Talagrand's Inequality] \label{lem:talagrand}
Let $c,r>0$ and $X$ be a non-negative random variable determined by $n$ independent trials satisfying the following conditions:
\begin{itemize}
\item ($c$-Lipschitz) changing the outcome of one trial can affect the value of $X$ by at most $c$; and 
\item ($r$-certifiable) for any $s>0$, if $X\geq s$ then there is a set of at most $rs$ trials whose outcomes certify that $X\geq s$.
\end{itemize}
Then for any $0\leq t\leq \lambda:=\bE[X]$, we have that
\[\bP\left[|X-\lambda|>t+60c\sqrt{r\lambda}\right]\leq 4e^{-\tfrac{t^2}{8c^2r\lambda}}.\]
\end{lem}

\section{Proof Overview} \label{sec:overview}

In the remainder of the paper (with the exclusion of the concluding remarks in Section~\ref{sec:conclude}), we fix  $1\leq j< k\in \NN$. For convenience, we also collect the following hierarchy of constants which is adhered to throughout; 
 \begin{linenomath}
\begin{equation}
   \label{eq:hierarchy} 
  1/n\ll \mu \ll \gamma  \ll 1/t \ll \beta \ll \eps, 1/k.
\end{equation}
\end{linenomath}
By this we mean that given any $\eps>0$ and $2\leq k \in \NN$, we can choose constants from right to left so that all required relations with these constants in our proof are satisfied. Indeed any statement concerning  these parameters in our proofs will be true under conditions of the form $c'\leq f(c_{1},\ldots,c_{\ell})$ where $c'$ is below $c_1,\ldots,c_{\ell}$ in the hierarchy. In order to highlight when we are appealing to this hierarchy, we will state that we are using that $c'\ll c_1,\ldots, c_\ell$. We will also fix  the following  key constants 
 \begin{linenomath}
\begin{equation} \label{eq:m}
T:=(t(k-1)+1), \hspace{4mm}     \tilde{m}:=t^{100k},  \hspace{4mm}    m:=T\cdot\tilde{m} \hspace{2mm}  \mbox{ and } \hspace{2mm} M:=T\cdot m=T^2\cdot \tilde{m}.
\end{equation}
\end{linenomath}
throughout. The definition of $\tilde{m}$ has some flexibility with the important thing being that it is much larger than $t$ but still bounded by some polynomial in $t$. The definitions of $T,m$ and $M$ are more rigid with it being crucial that $m/\tilde{m}=M/m=T$. Note that due to how the hierarchy is defined, we can assume that $\gamma \ll 1/M, 1/m, 1/\tilde{m},1/T$ and likewise with all constants to the left of $\gamma$ in the hierarchy. 
Finally, we will assume that $n\in (k-1)\NN$ and fix some $n$-vertex $k$-graph $G$ and a colouring $\chi:E(G)\rightarrow \NN$ satisfying the assumptions of Theorem~\ref{thm:main}.

\vspace{2mm}

Our proof will work through an application of the (lopsided) Lov\'{a}sz local lemma (Lemma~\ref{lem:llll}/ Corollary~\ref{L4}), which is a powerful tool for finding rainbow spanning structures in (hyper)graphs~\cite{dudek2012rainbow,albert1995multicoloured,erdos1991lopsided,coulson2019rainbow,ckpy2020rainbow,coulson2020rainbow}. Our probability space will be a choice of a uniformly random loose Hamilton cycle $H_*$ in $G$, noting that the existence of one such cycle is guaranteed by the fact that $G$ has minimum degree above the extremal threshold. Our bad events $\cA$ will encode the possibility of having two edges of the same colour in $H_*$ and thus our application of the local lemma will imply that with positive probability, the  random loose  Hamilton cycle $H_*$ is rainbow.

In order to appeal to the local lemma, we need to  bound from above probabilities of the form $\bP[A|\cap_{A'\in \cA'}\overline{A'}]$ for certain $\cA'\subset \cA$. Ignoring the conditioning of these events momentarily, we need an upper bound on the probability of some bad event $A$ occurring, where $A$ encodes the appearance of two edges of the same colour appearing in $H_*$. The key idea to estimate such  probabilities is to use \emph{switchings}. Broadly speaking, we do a double-counting argument of the pairs of loose Hamilton cycles $(H,H')$ of $G$ such that $H$ satisfies $A$, $H'$ does not satisfy it, and $H'$ can be obtained from $H$ by reshuffling a small number of disjoint subpaths. 
Exhibiting many such pairs for each $H$ in which the bad event $A$ occurs will allow us to conclude an upper bound on $\bP[A]$. 
The collection of subpaths that we will use to reshuffle the loose Hamilton cycle are what we will call a \emph{splitting} of the Hamilton cycle, as in the following definition.

\begin{dfn}[A Hamilton cycle splitting] \label{def:splitting}
    Given a loose Hamilton cycle $H$ in  $G$, an \emph{$m$-splitting} of $H$  is a collection of vertex-disjoint loose paths $\cP=\{P_1,\ldots, P_{m}\}$ in $G$ such that each path $P_i$ is a subgraph of $H$. We say the splitting is \emph{$t$-balanced} if each $P_i$ has length $t$ and we say that the splitting is \emph{$t$-bounded} if each $P_i$ has length \emph{at most} $t$.

    We define $V(\cP)=\cup\{V(P):P\in \cP\}$  to be vertices of $\cP$, $V^E(\cP):=\cup\{V^E(P):P\in \cP\}$ to be the \emph{endvertices} of  $\cP$ and $V^I(\cP):=V(\cP)\setminus V^E(\cP)=\cup\{V(P)\setminus V^E(P):P\in \cP\}$ to be the \emph{interior} vertices of $\cP$.
\end{dfn}

Given a splitting $\cP$ of a loose Hamilton cycle $H$ satisfying a bad event $A$, we will choose new edges of $G$ to piece back together a new loose Hamilton cycle $H'$ which avoids $A$. The edges that we use to do this will lie on the vertices $V(\cP)$ of the splitting and furthermore, they will be \emph{transverse}, that is, each new edge will contain at most one vertex from each of the paths in the splitting $\cP$. Using transverse edges in the new loose Hamilton cycle $H'$ will give us control over the colours that feature on new edges which, as we will see shortly, is essential to our approach. In general, we make the following definition of transverse sets and partitions. 

\begin{dfn}[Transverse with respect to a splitting] \label{def:transverse}
Given a loose Hamilton cycle $H$ in $G$ and an  $m$-splitting $\cP=\{P_1,\ldots, P_m\}$ 
of $H$, a set of vertices $S\subset V(\cP)$ is \emph{transverse} (with respect to $\cP$) if $|S\cap V(P)|\leq 1$ for all $P\in \cP$. A  partition of $V(\cP)$ is transverse (with respect to $\cP$) if each part of the partition is transverse.
\end{dfn}

Next we come to the notion of a \emph{switching}, which encapsulates how we reshuffle a loose Hamilton cycle $H$ on the vertices of some splitting $\cP$ to arrive at a new loose Hamilton cycle $H'$. The switching will depend on a path $P_0$ which is a subpath of the original loose Hamilton cycle $H$ that can be thought of as a ``troublesome" part of $H$. In applications, this path $P_0$ will contain (at least one of) the pair of monochromatic edges described by some bad event $A$ that occurs in $H$.

\begin{dfn}[Switchings] \label{def:switching}
Given a loose Hamilton cycle $H$ in $G$ and a loose Hamilton path $P_0$ of length $t$ with $P_0\subset H$, we define a  \emph{$(P_0,H)$-switching} of size $m$ to be a triple $(\cP,H',\cP')$ such that $\cP$ is a $t$-balanced  $m$-splitting of $H$ with $P_0\in \cP$, $H'$ is a loose Hamilton cycle of $G$, $\cP'$ is a $2t$-bounded $m$-splitting of $H'$, and the following conditions are satisfied. 
\begin{enumerate}[label=(\Roman*)]
\item \label{cond:fixed} $H\left[V(G)\setminus V^I(\cP)\right]=H'\left[V(G)\setminus V^I(\cP')\right]$;
\item \label{cond:splitP_0} $V(P_0)$ is transverse with respect to $\cP'$. 
\end{enumerate}

\end{dfn}

Condition \ref{cond:fixed} essentially says that when we switch from $H$ to $H'$, all of the altercation takes place on the vertices of a splitting $\cP$. Note also that condition~\ref{cond:fixed} necessarily implies that $V^E(\cP)=V^E(\cP')$. Condition \ref{cond:splitP_0} will imply that in the new loose Hamilton cycle $H'$, the vertices of $V(P_0)$ are far apart. In particular, this will imply that our bad event $A$ in $H$ will no longer occur in the new loose Hamilton cycle $H'$ as at least one of the pair of monochromatic edges described by $A$ will no longer feature in $H'$. 

As previously mentioned, we will use switchings to double count pairs $(H,H')$ of loose Hamilton cycles in $G$ in which some bad event $A$ occurs in $H$ but not in $H'$. This in turn will lead to upper bounds on the probability $\bP[A]$ of $A$ occurring. However, in order to appeal to the lopsided Lov\'asz local lemma (Lemma \ref{lem:llll}), what we actually want to bound is probabilities of the form $\bP[A|\cap_{A'\in \cA'}\overline{A'}]$ for any collection of events $\cA'$ that do not lie in the neighbourhood of $A$ in the dependency graph $\Gamma$ (see Definition \ref{def:dep}). As usual in the local lemma, in order to be effective, we will define $\Gamma$ \emph{locally} and so bad events $B$ that are adjacent to $A$ in $\Gamma$ will be those in which the edges involved in the event $B$ intersect the edges that are involved in the event $A$. In practice this will translate to having bad events $B$ that are adjacent to $A$ being those that involve edges intersecting the troublesome path $P_0$ and all other bad events might feature in the set $\cA'$. We will condition on some such collection of events $\cA'$ \emph{not} happening and in order to estimate $\bP[A|\cap_{A'\in \cA'}\overline{A'}]$, our double counting will focus in on loose Hamilton cycles  such that no event in $\cA'$ occurs. Thus we assume that no event in $\cA'$ occurs in $H$ and our switching will also need to produce some loose Hamilton cycle $H'$ in which no event in $\cA'$ occurs. This will be achieved by enforcing that the switching is \emph{feasible} as in the following definition.

\begin{dfn}[Feasible switchings] \label{def:feasible}
Given a loose Hamilton cycle $H$ in $G$ and a loose Hamilton path $P_0$ of length $t$ with $P_0\subset H$, we say that a $(P_0,H)$-switching $(\cP,H',\cP')$ is \emph{feasible} if $F':=H'[V(\cP')\setminus V(P_0)]$ is rainbow and no edge of $F'$ shares a colour with $H\left[V(G)\setminus V(\cP^I)\right]$. 
\end{dfn}

Note that in order for a switching to be feasible we do \emph{not} require that $H'$ is fully rainbow outside of $V(P_0)$. Indeed, it could be that there are edges of the same colour in $H'[V(G)\setminus V^I(\cP')]=H[V(G)\setminus V^I(\cP)]$. The key point is that in switching from $H$ to $H'$, we are not introducing any \emph{new} pairs of monochromatic edges, apart from possibly some that include vertices of $V(P_0)$. As we assume that no event in our collection $\cA'$ occurs in $H$, this will guarantee that none of the events in $\cA'$ will occur in $H'$ neither. 
There is the  exception of pairs of monochromatic edges that intersect $V(P_0)$ and might appear in $H'$ (and not in $H$) but, as alluded to above, we will define our dependency graph $\Gamma$ in such a way that these events will be adjacent to $A$ in $\Gamma$ and thus will not feature in $\cA'$. 

We are now in a position to state our first main lemma which precisely shows that if we can find many switchings for any loose Hamilton cycle $H$ and troublesome subpath $P_0$, then we will be able to prove the existence of a rainbow loose Hamilton cycle. The proof of Lemma \ref{lem:lll application}, which we give in Section \ref{sec:lll}, will follow from applying the local lemma (Corollary \ref{L4}) and the double counting argument sketched above. 

\begin{lem}[Many feasible switchings implies a rainbow Hamilton cycle] \label{lem:lll application}
Suppose that for  every loose Hamilton cycle $H$ in $G$ and loose path $P_0$ of length $t$ with $P_0\subset H$, there are at least $\gamma n^{m-1}$  feasible $(P_0,H)$-switchings of size $m$. Then $G$ contains a rainbow loose Hamilton cycle in $G$.
\end{lem}

With Lemma \ref{lem:lll application} in hand, it remains to find many feasible switchings for any loose Hamilton cycle $H$ and subpath $P_0$. To find a feasible switching, we first need to find a splitting $\cP$ of $H$ that contains $P_0$. The next definition captures the desired properties of such a splitting, that will allow us to use it in a feasible switching.

\begin{dfn}[Suitable splittings] \label{def:suitable}
 We say an $m$-splitting $\cP$  of a loose Hamilton cycle $H$  is \emph{suitable} for $P_0\in \cP$ if 
\begin{enumerate}[label=(\roman*)]
\item \label{suit:H cols} for any transverse $S\subset V(\cP)\setminus V(P_0)$ with $|S|=k-1$, there are at most $\eps m/4$ edges in $E_{G[V(\cP)]}(S)$ that share a colour with an edge in $H$; 
\item \label{suit:rb} 
for any transverse $X\subset V(\cP)\setminus V(P_0)$, there is no pair of edges $e,f\in E(G[X])$ with  $\chi(e)=\chi(f)$ and  $|e\cap f|= 1$; and 
\item \label{suit:disj} for any pair of \emph{disjoint} transverse edges $e,f\in E(G[V(\cP)\setminus V(P_0)])$ with $\chi(e)=\chi(f)$, we have that there are two distinct $P_i,P_{i'}\in 
\cP$ such that $P_x\cap e, P_x\cap f\neq \emptyset$ for $x=i,i'$.
\end{enumerate}
\end{dfn}

\begin{rem} \label{rem:suit}
 Note that condition~\ref{suit:disj} of Definition~\ref{def:suitable} in particular implies that for any transverse $X\subset V(\cP)\setminus V(P_0)$, there is no pair of disjoint edges $e,f\in E(G[X])$ with  $\chi(e)=\chi(f)$.
\end{rem}

Whilst the conditions for a splitting being suitable are somewhat technical, suffice it to say here that these conditions will guarantee that we can choose new  edges (all of which will be transverse with respect to the splitting $\cP$)  to contribute to the new loose Hamilton cycle $H'$ in our switching, in such a way that the switching will be feasible (Definition \ref{def:feasible}) and so not create any new unwanted pairs of monochromatic edges. 

On top of finding a suitable splitting, we also need to specify how we build our new $H'$ for our splitting and in particular, the endpoints of the new paths in the splitting $\cP'$ of $H'$. This will be given by a \emph{rerouting}.

\begin{dfn}[Rerouting] \label{def:reroute}
Given a loose Hamilton cycle $H$ in $G$ and an $m$-splitting $\cP=\{P_1,\ldots, P_m\}$ 
of $H$, a \emph{rerouting} of $\cP$ is a partition of $V^E(\cP)$ into pairs $\{\{v^i_a,v^i_b\}:i=1,\ldots, m\}$ such that connecting the paths in $H\left[V(G)\setminus V^I(\cP)\right]$ by identifying $v^i_a$ and $v^i_b$ for $i=1,\ldots,m$ results in a (single) loose cycle. 
\end{dfn}

Recall the definition of a partition  being transverse with respect to an $m$-splitting (Definition~\ref{def:transverse}) and note that if  we have a $t$-balanced $m$-splitting $\cP$ and $\cX=\{X_1,\ldots,X_{T}\}$  a  partition of $V(\cP)$ that is transverse with respect to $\cP$, then necessarily we have that each part in $\cX$ has size exactly
$m=T\cdot \tilde{m}$. Given a (suitable) splitting $\cP$ of a loose Hamilton cycle $H$, in order to find a feasible switching $(\cP,H',\cP')$ which uses $\cP$, we will first find transverse  partition $\cX$ of  $V(\cP)$ and then define our new loose Hamilton cycle $H'$ and corresponding splitting $\cP'$ in such a way that each path in $\cP'$ belongs to one part of the partition $\cX$. In particular, this will guarantee that each of the edges that features in $H'$ and not in $H$ will be transverse with respect to $\cP$. 
The next definition captures what properties we require from the partition of $V(\cP)$. 

\begin{dfn}[Viable partition] \label{def:viable}
For a $t$-balanced $m$-splitting $\cP$ of a loose Hamilton cycle $H$, we say a  partition $\cX=\{X_1,\ldots,X_{T}\}$ of $V(\cP)$  is \emph{viable} if it is transverse with respect to $\cP$  and 
\begin{enumerate}[label=(\alph*)]
\item \label{via:reroute} there exists a rerouting $\{\{v^i_a,v^i_b\}:i=1,\ldots, m\}$ of $\cP$ such that for each  $h\in [T]$ there are (exactly) $\tilde{m}$ indices $i\in [m]$ such that  $\{v^i_a,v^i_b\}\subset X_h$;
 and
\item \label{via:deg}$\delta_j(G[X_h])\geq (\delta_j^k(1)+\tfrac{\eps}{2})m^{k-j}$ for each $h\in [T]$. 
\end{enumerate}
\end{dfn}

Condition \ref{via:reroute} is necessary because we wish to build the paths in $\cP'$ so that each path is contained in one part of $\cX$ and so we certainly need a rerouting (which dictates the endpoints of the paths in $\cP'$) in which  pairs of vertices that feature in the rerouting lie in one part of the partition $\cX$. The condition that there are \emph{exactly} $\tilde{m}$ pairs in each part will imply that necessary divisibility conditions are adhered to so that we can cover the vertices of each part with disjoint paths between pairs of the rerouting. The condition \ref{via:deg} will be useful also to imply the existence of these disjoint paths. Indeed, without a minimum degree condition we will have no guarantee that $G[X_h]$ for $h\in[T]$ contains any sort of useful spanning structure for creating $H'$.

Our next major lemma, proved in Section \ref{sec:absorption}, shows that it in order to find a feasible switching, it is enough to find a suitable switching that admits a viable partition. 

\begin{lem}[Suitable splittings with viable partitions give feasible switchings] \label{lem:suitgivesfeas}
For a loose Hamilton cycle $H$ in $G$ and a path $P_0$ of length $t$ with $P_0\subset H$, if $\cP$ is a $t$-balanced $m$-splitting of $H$ that is suitable for $P_0$ and such that $V(\cP)$ admits a viable partition, then there is some (size $m$) feasible $(P_0,H)$-switching $(\cP,H',\cP')$.
\end{lem}

To prove this lemma, we will build the paths in $\cP'$ which define the new loose Hamilton cycle $H'$, by building paths in each part of the viable partition $\cX$, one at a time. Within each part $X_h$, at a high level, we use absorption techniques to find the relevant collection of disjoint paths that cover the vertices $X_h$. In fact, we are able to use a recent result of Alvarado, Kohayakawa, Lang, Mota and Stagni \cite{alvarado2023resilience}  (which is proven using absorption techniques) as a black box to find our desired paths. Still though, there is some technical work needed to get in to a position where we can apply their result, and we also need to prove that we can build our paths in such a way that the resulting switching is feasible, which we will be able to achieve by appealing to the conditions of  $\cP$ to be suitable.

Finally then, to complete the proof of Theorem \ref{thm:main}, it will remain to prove that we can indeed find suitable splittings that admit feasible partitions. This is the content of our last main lemma, Lemma \ref{lem:many suitable}, which is proven in Section \ref{sec:random sample} by sampling both a random splitting and a random partition and showing that with positive probability, all the required conditions will hold in these random samples. 

\begin{lem}[Many suitable splittings with viable partitions] \label{lem:many suitable}
For a loose Hamilton cycle $H$ in $G$ and a path $P_0$ of length $t$ with $P_0\subset H$, there exists $\gamma n^{m-1}$ distinct choices of $t$-balanced $m$-splittings $\cP$ of $H$ such that $\cP$ is suitable for $P_0$ and $V(\cP)$ admits a viable partition. 
\end{lem}

With these three main lemmas in hand, the proof of Theorem~\ref{thm:main} is immediate. We conclude this section by spelling it out in detail. 

\begin{proof}[Proof of Theorem~\ref{thm:main}]
By Lemma~\ref{lem:lll application}, it suffices to show that for any loose Hamilton cycle $H$ in $G$ and loose path $P_0\subset H$ of length $t$, there are at least $\gamma n^{m-1}$  feasible $(P_0,H)$-switchings of size $m$. Lemma~\ref{lem:suitgivesfeas} then implies that it is sufficient to find $\gamma n^{m-1}$ $t$-balanced $m$-splittings $\cP$ of $H$ with $P_0\in \cP$ and such that $\cP$ is suitable for $P_0$ and admits a viable partition. The existence of these splittings is then given by Lemma~\ref{lem:many suitable}. 
\end{proof}

\section{Many feasible switchings imply a rainbow Hamilton cycle} \label{sec:lll}

In this section we prove Lemma~\ref{lem:lll application}, which has a similar proof to that of~\cite[Lemma~2.2]{ckpy2020rainbow}.


\begin{proof}[Proof of \ref{lem:lll application}]
Let $H_*$ be a loose Hamilton cycle of $G$ chosen uniformly at random. Recall the hierarchy of the parameters: $1/n\ll\mu\ll \gamma\ll 1/m\ll 1/t$ in~\eqref{eq:hierarchy}.

For every $P_0$ loose path of length $t$ in $G$ and every $e_0,f_0\in E(P_0)$ with 
$|e_0\cap f_0|\leq 1$ and $\chi(e)=\chi(f)$ define the events $A(e_0,f_0;P_0)=\{P_0\subset H_*\}$. Similarly, for every pair of vertex-disjoint loose paths $P_0,Q_0$ of length $t$ in $G$ and every $e_0\in E(P_0)$ and $f_0\in E(Q_0)$ with $\chi(e)=\chi(f)$, define $B(e_0,f_0;P_0,Q_0)=\{P_0,Q_0\subseteq H_*\}$. 
Then we define 
 \begin{linenomath}
\begin{equation*}
    \begin{aligned}
        \mathcal{A}&=\{A(e_0,f_0;P_0): \bP[A(e_0,f_0;P_0)]>0\},\\
        \cB&={\{B(e_0,f_0;P_0,Q_0): \bP[B(e_0,f_0;P_0,Q_0)]>0\}},\\
        \cE&= \cA \cup \cB.
    \end{aligned}
\end{equation*}
\end{linenomath}
Note that if none of the events in $\cE$ hold, then $H_*$ is rainbow. Indeed, for every pair $e_0,f_0\in E(H_*)$ we have $|e_0\cap f_0|\leq 1$ and by avoiding the events in $\cE$, any such pair satisfies $\chi(e)\neq\chi(f)$. 

For each $A=A(e_0,f_0;P_0)\in \cA$, we write $V(A):=V(P_0)$ and for each $B=B(e_0,f_0;P_0,Q_0)$, we write $V(B):=V(P_0)\cup V(Q_0)$.
Consider the graph $\Gamma$ with $V(\Gamma)=\cE$ where for every $E_1,E_2\in \cE$ we have $E_1E_2\in E(\Gamma)$ if and only if $V(E_1)\cap V(E_2)\neq \emptyset$. We will choose values $p_\cA,p_\cB\in [0,1/2)$ 
such that $\Gamma$ is a $\textbf{p}$-dependency graph for $\cE$, 
where $p_A=p_\cA$ for all $A\in \cA$ and $p_B=p_\cB$ for all $B\in \cB$. For $\cX\in \{\cA,\cB\}$, define
 \begin{linenomath}
$$
d_{\cX} =\max_{E\in \cE} |\{X\in \cX: EX\in E(\Gamma)\}|.
$$
\end{linenomath}
To apply Corollary~\ref{L4}, it suffices to check that $p_{\cA}d_{\cA}+p_{\cB}d_{\cB} \leq 1/4$.

To bound the degrees in $\Gamma$, we bound the number of events whose vertex set intersects a given vertex $v\in V(G)$. 
Recall that  $T=t(k-1)+1$ is the number of vertices of a loose path of length $t$.

For every $v\in V(G)$, we claim that
 \begin{linenomath}
\begin{equation}\label{eq:A}
 |A\in \cA:v\in V(A)|\leq 2(k+1) T!\mu n^{T-1}.
\end{equation}
\end{linenomath}
Indeed, suppose that $A=A(e_0,f_0;P_0)\in \cA$ with $v\in V(A)$. There are at most $(k+1)\mu n^{2k-1}$ pairs of edges $e_0,f_0\in E(H_*)$ of the same colour with $|e_0\cap f_0|\leq 1$; indeed, there are at most $n^{k}$ choices for $e_0$, at most $(k+1)$ for $e_0\cap f_0$ and at most $\mu n^{k-1}$ for $f_0$, by~\ref{item:Global}. From them, at most $(k+1)\mu n^{2k-2}$ pairs satisfy $v\in e_0\cup f_0$ (as one can assume that $v\in e_0$ which gives at most $n^{k-1}$ choices for $e_0$). If $v\notin (e_0\cup f_0)$, then there are at most $T!n^{T-2k}$ ways to complete the path $P_0$ containing both edges, and if $v\in (e_0\cup f_0)$ there are at most $T!n^{T-2k+1}$ ways to do it. In total, $v$ intersects with at most $2(k+1)T!\mu n^{T-1}$  events in $\cA$, proving~\eqref{eq:A}.

For every $v\in V(G)$, we claim that
 \begin{linenomath}
\begin{equation}\label{eq:B}
 |B\in \cB:v\in V(B)|\leq 2 (T!)^2 \mu n^{2T-2}.
\end{equation}
 \end{linenomath}
    Indeed, suppose that $B=B(e_0,f_0;P_0,Q_0)\in \cB$ with $v\in V(B)$. If $v\in (e_0\cup f_0)$, we may assume that $v\in e_0$ and there are at most $n^{k-1}$ choices for $e_0$ and at most $\mu n^{k-1}$ for $f_0$ of the same colour, again by~\ref{item:Global}. Then there are at most $(T!)^2 n^{2T-2k}$ choices for $P_0$ and $Q_0$  containing $e_0$ and $f_0$ respectively. If $v\notin (e_0\cup f_0)$, then there are at most $\mu n^{2k-1}$ pairs $e_0,f_0$ of the same colour, and at most $(T!)^2 n^{2T-2k-1}$ choices for $P_0$ and $Q_0$  containing $e_0$ and $f_0$ respectively, such that $v\in P_0\cup Q_0$. In total, $v$ intersects at most  $2 (T!)^2 \mu n^{2T-2}$ events in $\cB$, proving~\eqref{eq:B}.

By the definition of $\Gamma$ and since $|V(E)|\leq 2T$ for all $E\in\cE$, there exists $C=C(k,t)$ such that
 \begin{linenomath}
\begin{align}\label{eq:deg}
d_{\cA} < C \mu n^{T-1} 
\quad \text{ and } \quad
d_{\cB} < C \mu n^{2T-2}. 
\end{align}
 \end{linenomath}

In order to bound $p_\cA$ and $p_\cB$ we will use switchings.
Recall that $p_\cA$ is an upper bound on $\bP[A\mid \cap_{E \in {\cE'}} \overline{E}]$
where $A \in \cA$ and 
 \begin{linenomath}
\begin{equation}\label{def:E'}
 {\cE'} \subseteq  \{E\in \cE\setminus \{A\} \,:\, AE \notin E(\Gamma)\}.
\end{equation}
\end{linenomath}

Fix $A=A(e_0,f_0;P_0)\in \cA$ and $\cE'$ as in~\eqref{def:E'}. Let $\cH=\cH(\cE')$ be the collection of  loose Hamilton cycles of $G$
that satisfy $\cap_{E \in {\cE'}} \overline{E}$, which we may assume  is non-empty.
Let $\cH_0 = \{H \in \cH: P_0 \subset H\}$. Consider the bipartite multigraph $\cF_A$
with parts $\cH_0$ and $\cH$ and
where an edge is added between $H \in \cH_0$ and $H'\in \cH$ for each feasible 
$(P_0,H)$-switching $Y=(\cP,H',\cP')$ of size $m$.

In fact, for every feasible $(P_0,H)$-switching $Y=(\cP,H',\cP')$, we have $H'\in \cH\setminus\cH_0$.
Indeed, by Definition~\ref{def:switching}, since $V(P_0)$ is transverse with respect to $\cP'$, $P_0$ is not in $H'$. Moreover, by Definition~\ref{def:feasible}, all edges in $H'$ that are not in $H$ either intersect $P_0$ or have a unique colour in $H'$. However, by the definition of $\Gamma$ and \eqref{def:E'}, all events that intersect $P_0$ are not in $\cE'$, thus $H'\in \cH$. 

Let $\delta_A$ be the minimum degree in  $\cF_A$ of the vertices in $\cH_0$ and $\Delta_A$ be the maximum degree in $\cF_A$ of the vertices in $\cH$. Double-counting the edges of $\cF_A$, we obtain
 \begin{linenomath}
\begin{equation}
    \bP[A \mid \cap_{E \in  {\cE'}}  {\overline{E}} ] = \frac{|\cH_0|}{|\cH|} \le \frac{\Delta_A}{\delta_A}.
\end{equation}
\end{linenomath}

It suffices to upper bound $\Delta_A$ and lower bound $\delta_A$.  The hypothesis of the lemma and the fact that every feasible switching produces a loose Hamilton cycle in $\cH$, imply that $\delta_A \ge \gamma n^{m-1}$.
To bound $\Delta_A$, we fix any $H'\in \cH$ 
and bound the number of pairs $(H,Y)$ where $H \in \cH_0$ 
and $Y=(\cP,H',\cP')$ is a feasible $(P_0,H)$-switching of size $m$. By the previous observation, it suffices to check it for $H'\in \cH\setminus \cH_0$.
By condition \ref{cond:splitP_0} in Definition~\ref{def:switching} , all vertices in $V(P_0)$ are in different parts of $\cP'$. We first construct $\cP'$. As $P_0$ is fixed by $A$, the number of ways to choose $T$ paths in $H'$ of length at most $2t$ and each one containing exactly one vertex in $V(P_0)$ is a constant only depending on $t,k$. Additionally, there are at most $n^{m-T}$ choices to complete $\cP'$ with $m-T$ paths in $H'$ (each factor of $n$ corresponding to a choice of starting point in $H'$ for one of the remaining paths in $\cP'$). Once $\cP'$ has been fixed, the number of ways to construct $\cP$ and $H$ is a constant only depending on $t,k,m$. So $\Delta_A \le C' n^{m-T}$, for $C'=C'(k,t,m)$.
We obtain
 \begin{linenomath}
\begin{align}\label{eq:pA}
\bP[A \vert \cap_{E \in  {\cE'}}  {\overline{E}} ] \leq C' \gamma^{-1} n^{1-T} =: p_\cA \;.
\end{align}
 \end{linenomath}

Similarly, we choose $p_\cB$ to be an upper bound on $\bP[B\mid \cap_{E \in  {\cE'}} \overline{E}]$
where $B=B(e_0,f_0;P_0,Q_0) \in \cB$ and ${\cE'}$ satisfying \eqref{def:E'} where $A$ is replaced by $B$. Define $\cH=\cH(\cE')$ as before and now let $\cH_{0} = \{H \in \cH: P_0,Q_0 \subset H\}$ and construct the bipartite multigraph $\cF_B$
with parts $\cH_{0}$ and $\cH$, 
where an edge is added between $H \in \cH_{0}$ and $H'\in \cH$ for each pair $(Y,Z)$ where 
\begin{itemize}
\setlength\itemsep{0.2cm}
    \item[-] $Y=(\cP,H_Y,\cP_Y)$ is a feasible 
$(P_0,H)$-switching of size $m$, such that $V(\cP)\cap V(Q_0)=\emptyset$, and
  \item[-] $Z=(\cQ,H',\cQ')$ is a feasible $(Q_0,H_Y)$-switching of size $m$, such that $V(\cQ)\cap V(P_0)=\emptyset$.
\end{itemize}
As in the previous observation, if there is an edge between $H$ and $H'$, then  $H'\in \cH\setminus \cH_{0}$. Define $\Delta_B$ and $\delta_B$ analogously to $\Delta_A$ and $\delta_A$ for the multigraph $\cF_B$. To compute $\delta_B$, note that some switchings are invalid as they violate the disjointness conditions. However, given $P_0$ and $Q_0$, there is some constant $\tilde{C}>0$ depending only on $t$ and $k$ such that there are at most  $\tilde{C}n^{m-2}$ many $(P_0,H)$-switchings $Y=(\cP,H_Y,\cP_Y)$ such that $V(\cP)\cap V(Q_0)\neq \emptyset$ and at most $\tilde{C}n^{m-2}$ many $(Q_0,H_Y)$-switchings $Z=(\cQ,H',\cQ')$ such that $V(\cQ)\cap V(P_0)\neq \emptyset$. 
Therefore a negligible proportion of the the pairs of switchings are disallowed and by the hypothesis of the lemma, $\delta_B\geq \tfrac{1}{2}\gamma^2 n^{2m-2}$.
Similarly as before, we have that $\Delta_B \le C'n^{2(m-T)}$. It follows that
 \begin{linenomath}
\begin{align}\label{eq:pB}
\bP[B \vert \cap_{E \in  {\cE'}}  {\overline{E}} ] \leq 2C' \gamma^{-2} n^{2-2T} =: p_\cB \;.
\end{align}
\end{linenomath}
Combining equation (\ref{eq:deg}), (\ref{eq:pA}) and (\ref{eq:pB}), and since $\mu\ll \gamma, 1/m, 1/t, 1/k$, we have that $p_{\cA}d_{\cA}+p_{\cB}d_{\cB} \leq 1/4$. The lemma now follows by applying Corollary~\ref{L4}.
    \end{proof}

\section{Suitable splittings with viable partitions give feasible switchings} \label{sec:absorption}

In this section, we will prove Lemma~\ref{lem:suitgivesfeas} showing that for a loose Hamilton cycle $H$ containing some length $t$ path $P_0$, if we have a $t$-balanced $m$-splitting $\cP$ that is suitable for $P_0$ and such that $V(\cP)$ admits a viable partition, then we can find some   feasible $(P_0,H)$-switching $(\cP,H',\cP')$. We therefore need to find a new loose Hamilton cycle $H'$ which is the same as $H$ outside of the vertex set $V^I(\cP)$. 
Recall that given the viable partition $\cX=\{X_1,\ldots, X_{T}\}$ of $V(\cP)$, there is some rerouting $\cR=\{\{v^i_a,v^i_b\}:i=1,\ldots,m\}$ of $\cP$  such that each $X_h$ contains exactly $\tilde{m}$ pairs in $\cR$. We will define $H'$ by 
considering each $X_h$ in the viable partition $\cX$ individually  
 and finding some collection $\cC_h$ of $\tilde{m}$ vertex-disjoint paths in $X_h$ that cover the vertex set of $X_h$ and such that each path $P\in \cC_h$ has as endpoints one of the pairs in $\cR$ that lies in $X_h$. Note that, due to the fact that $\cR$ is a rerouting (Definition~\ref{def:reroute}) and our requirement that the collections $\cC_h$ will cover all the vertices in $V(\cP)$, taking the paths in the collections along with the paths in $H[V(G)\setminus V^I(\cP)]$ will indeed define a new loose Hamilton cycle $H'$.  
 
 For each $h\in [T]$, in order to find the required path collection $\cC_h$ covering $X_h$, we will use the following lemma. 
 
 \begin{lem}[Path tilings] \label{lem:path tilings}
 Suppose that $G'\subset G$ is  a $k$-graph with $|V(G')|=m$ and $\delta_j(G')\geq (\delta_j^k(1)+\tfrac{\eps}{4})m^{k-j}$. Suppose also that there is some collection $\cR'=\{\{u^i_a,u^i_b\}:i\in[\tilde{m}]\}\subset \binom{V(G')}{2}$ of $\tilde{m}$ disjoint pairs of vertices in $G'$. Then there is a collection $\cC'=\{P_i:i\in[\tilde{m}]\}$ of $\tilde{m}$ vertex-disjoint (loose) paths in $G'$ each of length at most $2t$ and such that $V(G')=\cup_{i\in [\tilde{m}]} V(P_i)$ and $V^E(P_i)=\{u^i_a,u^i_b\}$ for each $i\in [\tilde{m}]$. Moreover, for any (2-uniform) graph $B$ with $V(B)=V(G')$ and $\Delta(B)\leq 2tk^2$, we can find a collection $\cC'=\{P_i:i\in[\tilde{m}]\}$ as above such that no edge of any path $P_i$ contains an edge of $B$. 
 \end{lem}
 
 We now show how to use Lemma~\ref{lem:path tilings} to establish Lemma~\ref{lem:suitgivesfeas}.

 \begin{proof}[Proof of Lemma~\ref{lem:suitgivesfeas}]
 We fix some path $P_0$ of length $t$ in $G$, some loose Hamilton cycle $H$ containing $P_0$ and some $t$-balanced $m$-splitting $\cP=\{P_0,\ldots,P_{m-1}\}$ of $H$, that is suitable for $P_0$ and such that $V(\cP)$ admits a viable partition 
 $\cX=\{X_1,\ldots,X_{T}\}$. By the definition of $\cX$ being viable (Definition~\ref{def:viable}), we have that there exists a rerouting $\cR=\{\{v^i_a,v^i_b\}:i=1,\ldots,m\}$ of $\cP$ and some partition of $\cR$ as $\cR=\cR_1\cup \cdots \cup \cR_{T}$ with $|\cR_h|=\tilde{m}$  and $\cR_h\subset \binom{X_h}{2}$ for each $h\in [T]$. Moreover, as $\cX$ is a viable (and hence transverse) partition of $V(\cP)$ of size $T$ we necessarily have that for each $h\in [T]$ and each $i=0,\ldots,m-1$, the intersection $X_h\cap V(P_i)$  contains a single vertex.  We label the vertices of each $X_h$ as $X_h=\{u_{0,h},\ldots,u_{m-1,h}\}$ so that $u_{i,h}\in V(P_i)$ for $i=0,\ldots, m-1$. We further set $\hat X_h:=X_h\setminus \{u_{0,h}\}=\{u_{i,h}:i\in[m-1]\}$ to be the vertex set obtained from $X_h$ by removing $u_{0,h}\in V(P_0)$.
 
 Our aim is to show that there is some feasible $(P_0,H)$-switching $(\cP,H',\cP')$. We will do this by finding a collection $\cC_h$ of  paths in $X_h$, 
  for $h\in [T]$ in succession. For each $h$, we will guarantee the following properties of $\cC_h$:
 \begin{enumerate}
 \item \label{item:path factor} $\cC_h$ is a \emph{path tiling} in $G[X_h]$, that is, the collection of paths in $\cC_h$ are vertex-disjoint and cover the vertex set $X_h$;
 \item \label{item:end pts} The set of pairs  of endpoints $\{V^E(P):P\in \cC_h\}$ is precisely set $\cR_h$ and each path in $\cC_h$ has length at most $2t$.
 \end{enumerate} 
 
For the next conditions we define $\hat{\cC_h}:=\{P[\hat X_h]:P\in \cC_h\}$ to be the collection of paths obtained from $\cC_h$ after deleting the vertex $u_{0,h}\in X_h$  and the edges that contain it. Note that by removing $u_{0,h}$, one of the paths in $\cC_h$ might split into two paths in $\hat\cC_h$. When we refer to edges in $\hat\cC_h$, we mean the collection of edges that feature on some subgraph belonging to the set $\hat\cC_h$. 
 \begin{enumerate}
\setcounter{enumi}{2}
  \item \label{item:no col outside} No edge in $\hat\cC_h$ shares a colour with an edge in $H[V(G)\setminus V(\cP^I)]$;
 \item \label{item: no col in prev} No edge in $\hat\cC_h$ shares a colour with an edge in  $\cup_{1\leq h'<h}\hat \cC_{h'}$;
 \item \label{item: internally rb} The subgraph of $G$ obtained by taking all edges in $\hat \cC_h$ is rainbow.
 \end{enumerate}
We claim that after finding such $\cC_h$, for all $h\in [T]$, we will be done. Indeed, note that the graph obtained by taking $H[V(G)\setminus V(\cP)]$ and all the paths in $\cP':=\cup_{h=1}^{T}\{P:P\in\cC_h\}$, defines a new loose Hamilton cycle $H'$ in $G$, using here the conditions \eqref{item:path factor} and \eqref{item:end pts} above, and the fact that $\cR$ is a rerouting as in Definition~\ref{def:reroute}. Moreover, by~\eqref{item:end pts}, note that $\cP'$ is an $m$-splitting of $H'$ which is $2t$-bounded. Note also that $H'[V(G)\setminus V^I(\cP')]$ is precisely the graph $H[V(G)\setminus V^I(\cP)]$ and that $V(P_0)$ is transverse with respect to $\cP'$ as the partition $\cX$ was transverse with respect to $\cP$ and all the paths in $\cP'$ lie in one part of $\cX$. Thus the triple $(\cP,H',\cP')$ defines a  $(P_0,H)$-switching of size $m$ (as defined in Definition~\ref{def:switching}). Furthermore, the graph $F':=H'[V(\cP')\setminus V(P_0)]$ is rainbow due to the conditions \eqref{item: no col in prev} and \eqref{item: internally rb} and no edge in $F'$ shares a colour with an edge in $H[V(G)\setminus V(\cP^I)]$ due to the conditions given by \eqref{item:no col outside}. Hence $(\cP,H',\cP')$ is a feasible $(P_0,H)$-switching, as required. 
 
It remains to prove that we can find the collections $\cC_h$ as claimed and as previously mentioned, we will show how to construct these sequentially for $h=1,\ldots,T$. So suppose that we are at the point of constructing $\cC_h$ for some $h\in [T]$ and that all $\cC_{h'}$ with $1\leq h'< h$ have already been constructed. Firstly, we let $G'\subset G[X_h]$ be the $k$-graph obtained from $G[X_h]$ by deleting any edge $e\subset \hat X_h=X_h\setminus \{u_{0,h}\}$ that shares a colour with some edge in $H$. Note that as $\cP$ is a suitable splitting of $H$ and $\cX$ is a viable partition, we have that $\delta_j(G')\geq (\delta_j^k(1)+\tfrac{\eps}{4})m^{k-j}$. Indeed, for any subset $S\subset X_h$ with $|S|=j$, property~\ref{via:deg} of Definition~\ref{def:viable} implies that $S$ is contained in at least $(\delta_j^k(1)+\tfrac{\eps}{2})m^{k-j}$ edges in $G[X_h]$ and property~\ref{suit:H cols} of Definition~\ref{def:suitable}  implies that at most $\tfrac{\eps}{4} m^{k-j}$ of these share a colour with some edge in $H$. We will find our paths in $G'$ and thus ensure that property \eqref{item:no col outside} holds. 

Next we  define a (2-uniform) graph $B_h$ with vertex set $V(B_h)=X_h=\{u_{0,h},\ldots,u_{m-1,h}\}$ and a pair of vertices $\{u_{i,h},u_{i',h}\}\in \binom{X_h}{2}$ being an edge of $B_h$ if and only if $0\notin \{i,i'\}$ and there is some $1\leq h'<h$ such that the pair $\{u_{i,h'},u_{i',h'}\}$ lies  in an  edge $e\in E(G[\hat X_{h'}])$ that features in $\hat\cC_{h'}$. We claim that $\Delta(B_h)\leq 2tk^2$. Indeed, for a vertex $u=u_{i,h}\in X_h$ we have that $\deg_{B_h}(u)=0$ if $i=0$ and otherwise we bound the number of choices for $i'$  with $\{u,u_{i',h}\}\in E(B_h)$ by a choice of $h'<h$ (at most $t(k-1)\leq tk$ choices), a choice of edge $e\in G[\hat X_{h'}]$ that lies in   $\hat \cC_h\subset \cC_h$ and  contains $u_{i,h'}$ (at most 2 choices) and a choice of vertex $u_{i',h'}$ that lies in $e$ (at most $k-1\leq k$ choices). 

By the `moreover' part of Lemma~\ref{lem:path tilings}, we can obtain a collection $\cC_h$ of $\tilde{m}$ paths  in $G'[X_h]$ that satisfy the conditions \eqref{item:path factor} and \eqref{item:end pts} above and such that no edge in a path in $\cC_h$ contains an edge of $B_h$. It remains to establish properties \eqref{item: no col in prev} and \eqref{item: internally rb}. For \eqref{item: internally rb}, note that as $\hat \cC_h$ is a collection of vertex-disjoint loose paths, any pair of edges $e,f\subseteq \hat X_h\subset V(\cP)\setminus V(P_0)$ that feature in $\hat \cC_h$ necessarily have that $|e\cap f|\leq 1$. Thus, as $\hat X_h$ is transverse with respect to $\cP$, we have that $\chi(e)\neq \chi(f)$ due to conditions~\ref{suit:rb} (for $|e\cap f|=1$) and~\ref{suit:disj} (for $|e\cap f|=0$, see Remark~\ref{rem:suit}) of Definition~\ref{def:suitable}  and the fact that $\cP$ is suitable for $P_0$. This establishes \eqref{item: internally rb}. Finally, for \eqref{item: no col in prev} suppose for a contradiction that there is some edge $e$ in $\hat\cC_h$  that shares a colour with some edge $f$  in $\hat \cC_{h'}$ for some  $1\leq h'<h$. As $\cP$ is suitable for $P_0$ and $e,f\subset V(\cP)\setminus V(P_0)$ are disjoint edges that are both transverse with respect to $\cP$ (as they lie in different parts of the transverse partition $\cX$), condition~\ref{suit:disj} of Definition~\ref{def:suitable} implies that there exists distinct $i,i'\in [m-1]$ such that $u_{i,h'},u_{i',h'}\subset f$ and $u_{i,h},u_{i',h}\subset e$. However, this implies that $\{u_{i,h},u_{i',h}\}\in E(B_h)$ contradicting the fact that no edge on a path in  $\cC_h$, contains an edge of $B_h$. This concludes the proof. 
 \end{proof}

\subsection{Path tilings} \label{sec:path tilings}
It remains to prove Lemma~\ref{lem:path tilings}. For this, we will use the following recent theorem due to Alvarado, Kohayakawa, Lang, Mota and Stagni \cite[Theorem 9.3]{alvarado2023resilience}  which shows that the $j$-degree threshold for a loose  Hamilton cycle coincides with the threshold for a loose Hamilton path between any pair of vertices. 

\begin{thm}[Hamilton paths] \label{thm:Ham path} For any $n_0\geq t/2$ such that $n_0-1\in (k-1)\NN$, if $G_0$ is an $n_0$-vertex $k$-graph with $\delta_j(G_0)\geq (\delta_j^k(1)+\tfrac{\eps}{8})n_0^{k-j}$, then for any pair of vertices $w^1,w^2\in V(G_0)$, there is a loose Hamilton path in $G_0$ with endpoints $w^1,w^2$.
\end{thm}

With Theorem~\ref{thm:Ham path} in hand, the simple idea behind proving Lemma~\ref{lem:path tilings} is to partition  $V(G')$ into $\tilde{m}$ sets $V_1,\ldots, V_{\tilde{m}}$ of  (roughly) equal sizes such that the collection of pairs $\cR'=\{\{u^i_a,u^i_b\}:i\in [\tilde{m}]\}$ is distributed amongst the vertex sets with $\{u^i_a,u^i_b\}\subset V_i$ for each $i\in[\tilde{m}]$. We can then hope to apply Theorem~\ref{thm:Ham path} to each $G'[V_i]$ individually to get a path $P_i$ in $G'[V_i]$ between $u^i_a$ and $u^i_b$ and thus obtain the desired path tiling $\cC'=\{P_i:i\in [{\tilde{m}}]\}$. In order to be able to apply Theorem~\ref{thm:Ham path} in each of the vertex sets $V_i$, we need the graphs $G'[V_i]$ to inherit the $j$-degree condition of $G'$. This is achieved by taking a random partition 
(apart from pairs in $\cR'$ which are placed in the desired set $V_i$). We also need the divisibility conditions that $|V_i|-1\in (k-1)\NN$ for each $i\in [\tilde m]$. When randomly distributing the vertices amongst the $V_i$, we will not be able to ensure such a divisibility condition but we overcome this by setting aside some vertex sets $W_i$ with $|W_i|=k-2$ for $i\in[\tilde m]$ which will be used to adjust the sets $V_i$ and establish the required divisibility conditions. In order to prove the `moreover' statement of Lemma~\ref{lem:path tilings}, we show that with high probability the random sets $V_i$ induce almost independent sets in $B$. If each $V_i$ was truly an independent set in $B$ then  the loose Hamilton paths output by Theorem~\ref{thm:Ham path} would not use hyperedges that contain edges of $B$. To turn the almost independent sets into truly independent sets in $B$, we will need to adjust our sets slightly by moving some limited number of vertices between parts of the partition. Moreover, by choosing our sets $W_i$ carefully to avoid edges of $B$, we can also maintain that the $B[V_i]$ are empty graphs after adding vertices from $W_i$ to amend the divisibility constraints. We now give the details.

\begin{proof}[Proof of Lemma~\ref{lem:path tilings}]
We will prove the whole statement including the `moreover' part and so we fix some $G'$ with vertex set $V=V(G')\subset V(G)$, some (2-uniform) graph $B$ with $V(B)=V$ and some collection of pairs of vertices $\cR'=\{\{u^i_a,u^i_b\}:i\in [\tilde{m}]\}$ as in the statement of the lemma. We let $R:=\cup_{i=1}^{\tilde{m}}\{u^i_a,u^i_b\}$ be the set of vertices that appear in a pair in $\cR'$. We will now define vertex sets $W_i,Y_i\subset V$ for $i=1,\ldots,\tilde{m}-1$ sequentially, and initiate the process by defining $W_0=\emptyset$ and  $Y_0=V\setminus R$. For each $i\in [\tilde{m}-1]$ we will define $W_i\subset Y_{i-1}$ to have $|W_i|=k-2$ and $Y_i$ will be defined as $Y_i=Y_{i-1}\setminus W_{i-1}$. This will guarantee that the sets $W_i$ are disjoint from each other and for each $i\in [\tilde{m}-1]$ we will have that $|Y_i|\geq m-2\tilde{m}-(k-2)(i-1)\geq (t-1)\tilde{m}$. 
We will also impose that $W_i\subset Y_{i-1}$ has the property that the only possible edges in $B[\{u^i_a,u^i_b\}\cup W_{i-1}\cup W_i]$ and $B[\{u^{i+1}_a,u^{i+1}_b\}\cup W_i]$ are  $\{u^i_a,u^i_b\}$ and $\{u^{i+1}_a,u^{i+1}_b\}$ respectively. Let us see that such a vertex set $W_i$ indeed exists. We start with $Y_{i-1}$ which has size at least $(t-1)\tilde{m}$ and remove all the vertices that are $B$-neighbours of a vertex in $\{u^i_a,u^i_b,u^{i+1}_a,u^{i+1}_b\}\cup W_{i-1}$ to get some vertex set $Y_{i-1}'\subset Y_{i-1}$ which has size at least $(t-1)\tilde{m}-(4+(k-2))2tk^2\geq (t-2)\tilde{m}$ as $ 1/t\ll1/k$ and $\tilde{m}=t^{100k}$ by \eqref{eq:hierarchy} and \eqref{eq:m}. As $\Delta(B)\leq 2tk^2$ there is some vertex set $Y''_{i-1}\subset Y'_{i-1}$ with $|Y''_{i-1}|\geq |Y'_{i-1}|/(2tk^2+1)\geq \tilde{m}/2tk^3\geq k-1$ such that $B[Y''_{i-1}]$ is empty. Choosing $W_i$ arbitrarily from this vertex set $Y''_{i-1}$ gives the desired properties, noting here that we have previously chosen $W_{i-1}$ so that the only possible edge in $B[\{u^i_a,u^i_b\}\cup W_{i-1}]$ is $\{u^i_a,u^i_b\}$. After choosing $W_i$, we define $Y_i=Y_{i-1}\setminus W_i$ by removing the vertices of $W_i$ from $Y_{i-1}$ and move onto the next step $i+1$ (or finish, if $i=\tilde{m}-1$). 

At the end of the process above we have defined $W_i$ for $i=0,1,\ldots,{\tilde{m}-1}$ and the vertex set $V':=Y_{\tilde{m}-1}$ which is the vertex set obtained from $V$ by removing  $R$ and all of the sets $W_i$. We further  define $W_{\tilde{m}}=\emptyset$. Note that $|V'|=T\cdot\tilde{m}-2\tilde{m}-(k-2)(\tilde{m}-1)=\tilde{m}(t-1)(k-1)+k-2$.  We proceed with the following claim, which will find a partition $V'=U_1\cup \cdots \cup U_{\tilde{m}}$  of the vertices $V'$ which have certain properties. For each $i\in [\tilde{m}]$ we will define $U_i^+:=U_i\cup \{u^i_a,u^i_b\}\cup W_i\cup W_{i-1}$.  Recall that by construction, there are no edges in $B[U_i^+\setminus U_i]$ with the possible exception of the edge $\{u^i_a,u^i_b\}$. We say that $U_i$ is \emph{good} if there is no edge in $B[U_i^+]$ different from $\{u^i_a,u^i_b\}$, \emph{bad} if there is at least one edge in $B[U_i^+]$ different from $\{u^i_a,u^i_b\}$. Moreover we say $U_i$ is \emph{very bad} if there is no choice of a single vertex $u\in U_i$ that is incident to all edges in $B[U^+_i]$ other than $\{u^i_a,u^i_b\}$. 
Recall that we have chosen $\beta$ such that $1/t \ll \beta \ll \eps$.

\begin{clm} \label{clm:partition Ui} There is a partition $V'=U_1\cup \cdots \cup U_{\tilde{m}}$  of the vertices $V'$ such that the following holds:
\begin{enumerate}[label=(\Alph*)]
\item \label{Ui size} $(t-1)(k-1)-\beta t\leq|U_i|\leq (t-1)(k-1)+\beta t$ for all $i\in [\tilde{m}]$;
\item \label{Ui degs} for any $S\in \binom{V'}{j}$ and $i\in [\tilde{m}]$, we have that $\delta_{G'}(S;U_i^+)\geq (\delta_j^k(1)+\tfrac{3\eps}{16})|U_i^+|^{k-j}$;  
\item \label{Ui v bad} none of the $U_i$ with $i\in[\tilde{m}]$ are very bad; and 
\item \label{Ui bad} at most $t^3k^3$ of the sets $U_i$ with $i\in [\tilde{m}]$ are bad.
\end{enumerate} 
\end{clm}
\begin{claimproof}
We will consider a random partition of the vertices of $V'$ and show that all the properties~\ref{Ui size},~\ref{Ui degs},~\ref{Ui v bad} and~\ref{Ui bad} hold with probability at least $4/5$. Thus with probability at least $1/5$ all four conditions will be satisfied  and so a partition as claimed does indeed exist. In more detail, each vertex $v\in V'$ is assigned to $U_i$ with probability $1/\tilde{m}$ for each $[\tilde{m}]$ and the allocations of different vertices is done independently. Note that the expected size of each $U_i$ is $\lambda:=(t-1)(k-1)+\tfrac{(k-2)}{\tilde{m}}\leq (t-1)(k-1)+1$. 

For part~\ref{Ui size}, note that for a fixed $i\in[\tilde{m}]$ we have that $|U_i|$ is binomially distributed and so by Chernoff bounds (Theorem~\ref{thm:chernoff}), 
\[\bP\left[|U_i|\notin[\lambda-\beta t-1,\lambda+\beta t]\right]\leq 2e^{-\frac{\beta^2t}{6k^2}}\leq e^{-\sqrt{t}}.\]
Thus the fact that with probability at least $9/10\geq 4/5$, ~\ref{Ui size} holds follows from a union bound over all $\tilde{m}=t^{100k}$ choices of $i\in [\tilde{m}]$, using here that $t$ is sufficiently large. 

For condition~\ref{Ui degs}, we fix  some $i\in [\tilde{m}]$  and some $j$-set $S\subset V=V(G')$. We define
$N(S):=\{A\subset V':A\cup S\in E(G')\}$ and note that $|N(S)|\geq (\delta^k_j(1)+\tfrac{\eps}{4}-\beta)m^{k-j}$. Indeed, we have that the number of edges in $G'$ containing $S$ is at least $(\delta^k_j(1)+\tfrac{\eps}{4})m^{k-j}$ due to the minimum degree condition on $G'$ and at most $|V\setminus V'|m^{k-j-1}\leq k\tilde{m}m^{k-j-1}\leq  \beta m^{k-j}$ of the sets $A\subset V(G')$ with $A\cup S\in E(G')$ are such that $A\nsubseteq V'$.  Now let $X_i(S)=|\{A\in N(S):A\subset U_i\}|$ be the random variable that counts the number of sets in $N(S)$ that end up in $U_i$. Each set in $N(S)$ is in $U_i$ with probability $(1/\tilde{m})^{k-j}$ and so $\lambda_i(S):=\bE[X_i(S)]\geq  (\delta^k_j(1)+\tfrac{\eps}{4}-\beta)(t(k-1))^{k-j}$. We will use Suen's inequality (Lemma~\ref{lem:suen}) to prove concentration for $\lambda_i(S)$ and so for each $A\in N(S)$ we let $I_A$ be indicator random variable for the event that $A\subset U_i$ and define a graph $\Gamma$ with vertex set $V(\Gamma)=N(S)$ such that $A\sim A'$ if and only if $A\cap A'\neq \emptyset$. Then $\Gamma$ is a strong dependency graph for the family of indicator random variables $\{I_A\}_{A\in N(S)}$ and  following the notation of Lemma~\ref{lem:suen}, we have that 
\[\Delta_{X_i(S)}= \sum_{\ell=1}^{k-j-1} \sum_{\substack{(A,A')\in N(S)^2 \\ |A\cap A'|= \ell}}\bE[I_{A}I_{A'}]\leq \sum_{\ell=1}^{k-j-1}\binom{2(k-j)-\ell}{\ell}m^{2(k-j)-\ell}\left(\frac{1}{\tilde{m}}\right)^{2(k-j)-\ell}\leq (2k)^{3k}t^{2(k-j)-1},\]
and similarly \[\delta_{X_i(S)}=\max_{A\in N(S)}\sum_{A':|A'\cap A|\geq 1}\bE[I_{A'}]\leq \sum_{\ell=1}^{k-j-1}\binom{k-j}{\ell} m^{k-j-\ell}\left(\frac{1}{\tilde{m}}\right)^{k-j}\leq (2tk)^{3k}m^{-1}.\]
By Suen's inequality (Lemma~\ref{lem:suen}), we have that $\bP[X_i(S)\leq \lambda_i(S)-\beta (t(k-1))^{k-j}]\leq e^{-\sqrt{t}}$ and so a union bound over all $m^{j}=(t^{100k}T)^j$ choices of $S\subseteq V$ and all $\tilde{m}=t^{100k}$ choices of $U_i$ gives that with probability at least $9/10$ we have that $X_i(S)\geq \lambda_i(S)-\beta (t(k-1))^{k-j}$ for each $i\in [\tilde{m}]$ and $S\in \binom{V}{j}$, using here that $t$ is sufficiently large. Moreover, we showed that condition~\ref{Ui size} holds with probability at least $9/10$ and so with probability at least $4/5$, both~\ref{Ui size} and the lower bounds on all $X_i(S)$ hold. This establishes the degree conditions~\ref{Ui degs} by using our lower bound on the values $\lambda_i(S)$, the upper bound on the sizes of $|U_i|$ (and the resulting upper bounds on $|U_i^+|$) from~\ref{Ui size} and the fact that $\beta\ll \eps$ in our hierarchy \eqref{eq:hierarchy}.

It remains to establish conditions~\ref{Ui v bad} and~\ref{Ui bad}. Again, note that we defined the sets $W_i$ in such a way that the only possible edge in $B[U_i^+\setminus U_i]$ is $\{u^i_a,u^i_b\}$  for all $i\in[\tilde{m}]$. 
For each $i\in [\tilde{m}]$, it suffices to consider edges in $B[U_i]$  and edges in $B$ between $U_i$ and $U_i^+\setminus U_i$. Fixing some $i\in [\tilde{m}]$ and considering edges in $B[U_i]$, we note that by the maximum degree condition on $B$, we have that $B[V']$ contains at most $2tk^2m$ edges and each edge lands in  $U_i$ with probability $1/\tilde{m}^2$. Similarly for edges between $U_i^+\setminus U_i$ and $U_i$, we have that there at at most $(2k-2)2tk^2$ edges between $U_i^+\setminus U_i=\{u^i_a,u^i_b\}\cup W_i\cup W_{i-1}$ and $V'$ and each such edge lies in $B[U_i^+]$ with probability $1/\tilde{m}$. Therefore the expected number of edges  with at least one endpoint in $U_i$ is at most $2tk^2m\cdot(1/\tilde{m})^2+(2k-2)2tk^2\cdot(1/\tilde{m})\leq 6t^2k^3/\tilde{m}$ and by Markov's inequality, we have that $\bP[U_i \mbox{ is bad}]\leq 6t^2k^3/\tilde{m}$ for each $i\in [\tilde{m}]$. Similarly, the probability that $U_i$ is very bad can be upper bounded by the expected number of pairs of edges in $B[U_i^+]$, each having an endpoint in $U_i$ using Markov's inequality. We get that 
\begin{linenomath}
\begin{align*}
\bP[U_i \mbox{ is very bad}]\leq \frac{(2tk^2m)^2}{
\tilde{m}^4}+ \frac{2tk^2m\cdot(2k-2)2tk^2}{
\tilde{m}^3} + \frac{((2k-2)2tk^2)^2}{
\tilde{m}^2}
\leq \frac{40k^6t^4}{\tilde{m}^2} ,
\end{align*}
 \end{linenomath}
 where the first summand accounts for vertex-disjoint pairs of edges in $B[U_i]$,  the second one accounts for an edge in $B[U_i]$ paired with a disjoint edge between $U_i$ and $U_i^+\setminus U_i$ and the last one counts pairs of edges which are both between $U_i$ and $U_i^+\setminus U_i$, and which by definition of very bad, do not intersect in $U_i$. Note that one of these three scenarios must necessarily occur for $U_i$ to be very bad as otherwise there would be one vertex in $U_i$ incident to all edges in $B[U_i^+]$ other than $\{u_i^a,u_i^b\}$. 
 Thus, by linearity of expectation, the expected number of very bad sets $U_i$ is less than $40k^6t^4/\tilde{m}<1/5$ and the fact that~\ref{Ui v bad} holds with probability at least $4/5$ follows from Markov's inequality. Similarly, the expected number of bad $U_i$ is at most $6t^2k^3$ by linearity of expectation, and so $\bP[\neg~\ref{Ui bad}]\leq 6/t\leq 1/5$ as required. 
\end{claimproof}

\vspace{2mm}

We now fix some partition $V'=U_1\cup \cdots\cup U_{\tilde{m}}$ as in Claim~\ref{clm:partition Ui} and we will first slightly adjust the partition so that all the sets become good. By condition~\ref{Ui bad} of Claim~\ref{clm:partition Ui} there is some set of indices $I_B\subset [\tilde{m}]$ with $|I_B|\leq t^3k^3$ and such that $U_{i'}$ is good for all $i'\in [\tilde{m}]\setminus I_B$. For each bad index $i\in I_B$ as $U_i$ is not very bad due to~\ref{Ui v bad}, we have that there is some vertex $z_{i}\in U_{i}$ such that $U'_{i}:=U_{i}\setminus \{z_{i}\}$ is \emph{good}, that is, such that the only possible edge in $B[U_i'\cup \{u^i_a,u^i_b\}\cup W_{i}\cup W_{i-1}]$ is $\{u^i_a,u^i_b\}$. Moreover, there are at least $\tilde{m}-t^3k^3-4tk^2\geq \tilde{m}/2$ indices $i'\in [\tilde{m}]\setminus I_B$ such that adding $z_i$ to $U_{i'}$ results in a set $U'_{i'}=U_{i'}\cup\{z_i\}$ which is good. Indeed, $z_i$ has at most $2tk^2$ neighbours in $B$ and each  neighbour of $z_i$ rules out at most $2$ possible choices for such an $i'$ (the 2 here is necessary in case the neighbour lies in some set $W_{i''}$). 
Therefore for each $i\in I_B$ corresponding to a bad set $U_i$ we can find a $i'=i'(i)$ such that adding $z_i$ to $U_{i'}$ results in both $U'_i:=U_i\setminus \{z_i\}$ and $U'_{i'}:=U_{i'}\cup \{z_i\}$ being good and we can do this for all $i\in I_B$ in such a way that $i'(i_1)\neq i'(i_2)$ for distinct $i_1,i_2\in I_B$. 
Let $V'=U'_1\cup \ldots \cup U'_{\tilde{m}}$ be the resulting partition after switching all the $z_i$. Here, for any set $i_0\in [\tilde{m}]\setminus I_B$ such that ${i_0}$ is not chosen as $i'(i)$ for some $i\in I_B$, we simply leave $U_{i_0}$ unchanged and define $U'_{i_0}:=U_{i_0}$.

Next we will distribute vertices in the sets $W_i$ to get a partition of $V=V(G')$ as $V_1\cup \cdots \cup V_{\tilde{m}}$ with $U_i'\cup \{u^i_a,u^i_b\}\subset V_i$  and $|V_i|-1\in (k-1)\NN$ for each $i\in [\tilde{m}]$. We define $V_i$ for $i=1,\ldots,\tilde{m}$ sequentially via the following process, initiating with $W^-_0:=W_0=\emptyset$ and $V_0=\emptyset$. Suppose we are at step $i\in[\tilde{m}]$ and that $V_{i'}$ and $W^-_{i'}$ are defined for all $0\leq i'<i$. Then we let $k_i:=(|U'_i|+|W^-_{i-1}|+1) \mod (k-1)$ and note that $k_i\in \{0,\ldots,k-2\}$. We define $W_i^+\subset W_i$ to be an empty set if $k_i=0$ and some arbitrary set of size $k-1-k_i$ if $k_i\neq 0$. Furthermore, we define $V_i$ to be $V_i:=U_i'\cup W_{i-1}^-\cup \{u^i_a,u^i_b\}\cup W_i^+$ and $W_i^-:W_i\setminus W_i^+$. We then continue to step $i+1$ or finish if $i=\tilde{m}$. Note that, for all $i\in [\tilde{m}-1]$ this process will succeed at step $i$ and result in a set $V_i$ such that $|V_i|\mod (k-1)=|U'_i|+|W^-_{i-1}|+2+(k-1-k_i) \mod (k-1)$ and so  we indeed have that $|V_i|-1\in (k-1)\NN$. For $i=\tilde{m}$ we claim that $k_i=0$ and so the algorithm also succeeds and $|V_{\tilde{m}}|-1\in (k-1)\NN$ (if $k_{\tilde{m}}\neq 0$ then the algorithm would stall as $W_{\tilde{m}}=\emptyset$ and so an appropriate $W_{\tilde{m}}^+$ would not exist). To see it, observe  that $U_{\tilde{m}}'\cup W_{\tilde{m}-1}^-\cup\{u^{\tilde{m}}_a,u^{\tilde{m}}_b\}=V\setminus (\cup_{i=1}^{\tilde{m}-1}V_i)$ and for each $i\in [\tilde{m}-1]$ we have that $|V_i|=1 \mod (k-1)$ from above. Thus as $|V|=T\cdot\tilde{m}= (t(k-1)+1)\tilde{m}=\tilde{m} \mod (k-1)$, we get that $|U_{\tilde{m}}'|+|W_{\tilde{m}-1}^-|+2=\tilde{m}-(\tilde{m}-1)=1\mod (k-1)$ and indeed $k_{\tilde{m}}=0$ as required.

Note that for each $i\in [\tilde{m}]$, the vertex set $V_i$ and the set $U_i^+$ defined via Claim~\ref{clm:partition Ui} differ by at most $2(k-1)+1$ vertices, and so condition~\ref{Ui degs} from Claim~\ref{clm:partition Ui} implies that $\delta_j(G'[V_i])\geq  (\delta_j^k(1)+\tfrac{\eps}{8})|V_i|^{k-j}$. Note also that $|V_i|\geq |U_i|\geq (t-1)(k-1)-\beta t\geq t/2$ by part~\ref{Ui size} of Claim~\ref{clm:partition Ui}. We are therefore in a position to apply Theorem~\ref{thm:Ham path} in $G'[V_i]$ to obtain some  loose Hamilton path $P_i$ in $G'[V_i]$ between $u^i_a$ and $u^i_b$. Doing this for each $i\in [\tilde{m}]$ results in some collection $\cC':=\{P_i:i\in [\tilde{m}]\}$ which we claim satisfies the required conditions. Indeed, part~\ref{Ui size} of Claim~\ref{clm:partition Ui} and how we formed $V_i$ implies that we certainly have $|V_i|\leq tk$ and so $P_i$ will have length less than $2t$, for each $i\in[\tilde{m}]$. Moreover, for each such $i\in[\tilde{m}]$, due to how we defined $V_i$ (and $U_i'$), we have that the only possible edge in $B[V_i]$ is $\{u^i_a,u^i_b\}$ and the pair $\{u^i_a,u^i_b\}$ does not lie in any edge of $P_i$ as they are the endpoints of the path. This concludes the proof.
\end{proof}

\section{Many suitable splittings with viable partitions}  \label{sec:random sample}

In this section, we prove Lemma~\ref{lem:many suitable} showing that for a fixed loose Hamilton cycle $H$ and loose path $P_0\subset H$, there are many suitable splittings containing $P_0$ that admit viable partitions. We encourage the reader to remind themselves of the definitions of suitable (Definition~\ref{def:suitable}) and viable (Definition~\ref{def:viable}) in order to follow the contents of this section.  Lemma~\ref{lem:many suitable} will be achieved by random sampling, showing that with probability bounded away from $0$, a randomly sampled $m$-splitting will suffice. As there are many $m$-splittings containing $P_0$, indeed there are $\Omega(n^{m-1})$, we thus obtain the desired lower bound. In fact, instead of sampling an $m$-splitting, we will use a slightly more general model to sample, with more independence which will be useful in the analysis. Moreover, we will achieve our aim with two rounds of randomness, first showing that our random sampling results in a suitable $m$-splitting with probability bounded from $0$ and then showing that such an $m$-splitting admits a viable partition which will be found again by appealing to a random choice of partition. In order to set ourselves up for this second round of randomness, we need slightly more from our random $m$-splitting than it just being suitable, as the following lemma details (see \eqref{eq:m} for the definition of $M$).

\begin{lem} \label{lem:random split}
 Let $G$ be as described in Theorem~\ref{thm:main}.
For a loose Hamilton cycle $H$ in $G$ and a path $P_0\subset H$ of length $t$, let $E_p\subset E(H)$ be a random subset of $E(H)$ obtained by sampling each edge independently with probability $p:=\tfrac{(m-1)(k-1)}{n}$ and $\cP_p$ be the collection of length $t$ increasing paths in $H$  starting from an edge in $E_p$. Then with probability at least $1/m$, we have that $P_0\notin \cP_p$, the collection of paths $\cP_*:=\cP_p\cup\{P_0\}$ is a $t$-balanced $m$-splitting\footnote{In particular the paths in $\cP_p$ are disjoint from $P_0$.} of $H$ suitable for $P_0$ and $\delta_j(G[V(\cP_*)])\geq (\delta_j^k(1)+\tfrac{3\eps}{4})M^{k-j}$.
\end{lem}

\begin{proof}
Before embarking on the proof, we make some useful definitions and observations. We say that a pair of vertices $u,v\in V(G)$ are \emph{close} if there are edges $e,f\in E(H)$ with $u\in e, v\in f$ and such that there is a path of length at most $2t+1$ in $H$ containing both $e$ and $f$. We say that $u,v\in V(G)$ are \emph{spread} if they are not close and we say that a set $S\subset V(G)$ is spread if there are no pairs of close vertices contained in $S$. We say that a subset $S\subset V(G)$ is \emph{almost spread} if there is \emph{exactly one} pair of vertices $u,v\in S$ that are close  but all other pairs of vertices in $S$ are spread. 

Now note that for a vertex $v\in V(G)$, the probability that $v$ belongs to the set $V(\cP_p)$ is  $1-(1-p)^{t+1}$ if $v$ is a vertex of degree $2$ in $H$ and $1-(1-p)^t$ if $v$ is not a  vertex of degree $1$ in $H$. In either case we have that 
 \begin{linenomath}
\begin{equation} \label{eq:single vx prob}
\frac{(1-\beta)M}{n}\leq \bP[v\in V(\cP_p)]\leq \frac{2tmk}{n}, 
\end{equation}
\end{linenomath}
recalling that $\beta$ is fixed so that $ 1/t\ll \beta\ll \eps$ in our hierarchy \eqref{eq:hierarchy} and $m$ is much larger than $t$; see \eqref{eq:m}.  
For a spread set of vertices $S\subset V(G)$ note that the events that each $v\in V(\cP_p)$ are independent for $v\in S$ and so $\bP[S\subset V(\cP_p)]=\prod_{v\in S}\bP[v\in V(\cP_p)]$.

We collect the following events:
\begin{linenomath}
\begin{align*}
\cE_1&:=\{V(\cP_p) \mbox{ contains a spread } S \mbox{ with }|S|=k-1 \mbox{ and there are at least }\eps m/4 \mbox{ edges  in }E_{G[V(\cP_*)]}(S)\\ & \qquad\mbox{that share a colour with an edge in }H\}, \\
\cE_2&:= \{V(\cP_p) \mbox{ contains a pair } e,f\in E(G) \mbox{ with }\chi(e)=\chi(f), |e\cap f|\leq 1 \mbox{ and with }e\cup f \mbox{  spread}\}, \\
\cE_3&:= \{V(\cP_p) \mbox{ contains a pair } e,f\in E(G) \mbox{ with }\chi(e)=\chi(f), e\cap f=\emptyset \mbox{ and }e\cup f \mbox{ almost spread}\}, \\
\cE_4&:= \{\delta_j(G[V(\cP_*)])< (\delta_j^k(1)+\tfrac{3\eps}{4})M^{k-j}\}, \mbox{ and } \\
\cE_5&:=\{V(\cP_*) \mbox{ contains a pair of close vertices in distinct paths in }\cP_*\}.
\end{align*}
\end{linenomath}
The proof of the lemma will be deduced from the following claim. 

\begin{clm} \label{clm:bad events}
For each $i\in [5]$, we have that $\bP[\cE_i]\leq 1/m$. 
\end{clm}

To see that Claim~\ref{clm:bad events} implies the lemma, first note that  by Lemma~\ref{lem:hit expected} with probability at least $1/(4\sqrt{m})\geq 6/m$, we have that $|E_p|=\bE[|E_p|]=m-1$ and hence $|\cP_*|=m$.
Thus, with probability at least $6/m-5\cdot (1/m)=1/m$, we have that $|\cP_*|=m$ and none of the events $\cE_1,\ldots, \cE_5$ hold, by Claim~\ref{clm:bad events}. In particular, event $\cE_5$ not holding implies that there are no pairs of  intersecting paths in $\cP_*$ and so $\cP_*$ is indeed a $t$-balanced $m$-splitting of $H$. The degree condition in the statement of the lemma is given directly by event $\cE_4$ not holding and it remains to establish that $\cP_*$ is suitable for $P_0$. For this, note that any transverse $Y\subset V(\cP)\setminus V(P_0)$ is necessarily spread by event $\cE_5$ not happening and so condition~\ref{suit:H cols} of being suitable (Definition~\ref{def:suitable}) follows from the absence of event $\cE_1$, whilst condition~\ref{suit:rb} follows from the absence of event $\cE_2$. Finally, for condition~\ref{suit:disj}  consider two  disjoint edges $e,f\in G[V(\cP_*)\setminus V(P_0)]$ with $\chi(e)=\chi(f)$ and suppose that there are \emph{not} two distinct $P_i, P_{i'}\in \cP_*\setminus P_0$ with $P_x\cap e, P_x\cap f\neq \emptyset$ for $x=i,i'$.  Then $e\cup f$ intersects at least $2k-1$ members of $\cP_*$ and thus $e\cup f$ intersects at most one $P_i\in \cP_*$ in two vertices (and no $P_i$ in more than 2 vertices).  Due to the absence of event $\cE_5$, there is no pair of close vertices in distinct paths in $\cP_*$. Hence we have that  $e\cup f$ can contain at most one pair of close vertices (coming from the at most one path that it intersects in two vertices), that is,  we must have that $e\cup f$ is spread or almost spread. If the former, then we get a contradiction due to $\cE_2$ not happening whilst if the latter is the case, then we contradict that $\cE_3$ did not happen. Thus~\ref{suit:disj} indeed occurs and $\cP_*$ is suitable for $P_0$, as required.


\vspace{2mm}

It remains to establish Claim~\ref{clm:bad events} which we now do by considering each event separately. 

\subsection*{Bounding $\bP[\cE_1]$.}
Let $\chi(H)$ denote the colours featuring on edges in $H$ and $E_{\chi(H)}:=\{e\in E(G): \chi(e)\in \chi(H)\}$ be the set of edges that share a colour with an edge of $H$. 
  Let us fix $\alpha:=\eps/(32tk)$ and say that a $(k-1)$-set $S\in \binom{V(G)}{k-1}$ of vertices in $G$ is \emph{bad} if $|E_G(S)\cap E_{\chi(H)}|\geq \alpha n$ and \emph{good} otherwise. We begin by showing that  the probability that there exists some  bad spread set  $S\in \binom{V(G)}{k-1}$ such that  $S\subset V(\cP_p)$ is at most $1/2m$. Indeed,  we have that $|E_{\chi(H)}|\leq \frac{n}{k-1}\cdot \mu n^{k-1}\leq \mu n^k$ using that  $H$ has $n/(k-1)$ edges (and so at most this many colours featuring in $\chi(H)$) and the global bound~\ref{item:Global} of Theorem~\ref{thm:main}. Hence the number of  bad spread $(k-1)$-sets in $G$ is at most $(\mu k n^{k-1})/\alpha$. Each such set $S$  appears in $V(\cP_p)$ with probability at most $(2tmk/n)^{k-1}$ using the upper bound in \eqref{eq:single vx prob} and the fact that $S$ is spread and so each vertex of $S$ appears independently. Hence, by Markov's inequality the probability that there exists a bad spread set in $V(\cP_p)$ is at most 
  \[\frac{\mu k n^{k-1}}{\alpha}\cdot \frac{(2tmk)^{k-1}}{n^{k-1}}\leq \frac{1}{2m},\]
  using that $\mu\ll 1/m,  1/t,1/k, \eps$ and $\alpha=\eps/(32tk)$. 
  
 Now for a good $(k-1)$-set $S\in \binom{V(G)}{k-1}$, let $Y(S):=|E_{G[V(\cP_*)]}(S)\cap E_{\chi(H)}|$ be the random variable which counts the number of edges in $G[V(\cP_*)]$ containing $S$ that share a colour with an edge of $H$. We will show that the probability that there is some good spread $S$ that appears in $V(\cP_p)$ with $Y(S)\geq \eps m/4$ is at most $1/(2m)$. Hence with probability at least $1-1/m$ all the spread $(k-1)$-sets in $V(\cP_p)$ are good and have a bounded $Y(S)$, establishing Claim~\ref{clm:bad events} for $\cE_1$. 
 
 In order to measure the probability of some good $S$ having large $Y(S)$, consider the subset $N(S)$ obtained from  $\{e\setminus S : e\in E_G(S)\cap E_{\chi(H)}\}\subset V(G)$ by removing all vertices that are close to a vertex in $S\cup V(P_0)$. Then consider the random variable $Z(S):=|N(S)\cap V(\cP_p)|$. As $Y(S)\leq Z(S)+ (2t+1)\cdot 2k^2 \leq Z(S)+\eps m/8=:Z'(S)$ for any good $S$, we will instead bound the probability that the random variable $Z'(S)$ is larger than $\eps m/4$. For bounding $Z'(S)$, note that for any good spread $S$, we have that $Z'(S)$ is $tk$-Lipschitz and $1$-certifiable. Moreover, using \eqref{eq:single vx prob} and the definition of $\alpha$ we have that $\lambda:=\bE[Z'(S)]\leq \alpha n \cdot (2tmk)/n +\eps m/8\leq  3\eps m/16$. Hence by Talagrand's Inequality (Lemma~\ref{lem:talagrand}), 
 \[\bP\left[Z'(S)\geq \eps m/4\right]\leq \bP\left[Z'(S)\geq \eps m/8 +60tk\sqrt{\lambda}\right]\leq 4e^{-\tfrac{\eps^2m^2}{512t^2k^2\lambda}}\leq m^{-(k+2)},\]
 as $1/m\ll 1/t,1/k,\eps$; see \eqref{eq:hierarchy}. Moreover, due to the fact that we removed vertices close to $S$ in order to get $N(S)$, the events $\{S\subset V(\cP_p)\}$ and $\{Z'(S)\geq \eps m/4\}$ are independent and so occur concurrently with probability at most $(2tmk/n)^{k-1}\cdot (1/m)^{k+2}$.  The required bound on the probability of some good spread $S\subset V(\cP_p)$ having $Y(S)\geq \eps m/4$ thus follows from a union bound over the at most $n^{k-1}$ choices for good spread $S$ in $V(G)$.

\subsection*{Bounding $\bP[\cE_2]$.}  For a pair of edges $e,f\in E(G)$ with $e\cup f$ being spread and $|e\cap f|\leq 1$ we have that $\bP[e\cup f\subset V(\cP_p)]\leq (2tmk/n)^{2k-1}$ by \eqref{eq:single vx prob}. The number of such pairs $e,f$ with $\chi(e)=\chi(f)$ is at most $\mu n^{2k-1}$ due to the  bound~\ref{item:Global} on the colouring in Theorem~\ref{thm:main}. Hence, by a union bound, 
$\bP[\cE_2]\leq \mu n^{2k-1}\cdot \left(\frac{2tmk}{n}\right)^{2k-1}\leq \frac{1}{m},$
using here that $\mu\ll 1/m,1/t, 1/k$. 

\subsection*{Bounding $\bP[\cE_3]$.} The number of pairs of disjoint edges $e,f\in E(G)$ with $\chi(e)=\chi(f)$ is at most  $\mu n^{2k-1}$ due to the global bound~\ref{item:Global} on the colouring in Theorem~\ref{thm:main}. If such a pair is almost spread, then $\bP[e\cup f\subset V(\cP_p)]\leq (2tmk/n)^{2k-1}$. Indeed this follows from \eqref{eq:single vx prob} and the fact that we can remove one vertex from $e\cup f$ to get a spread set. Hence by a union bound, $\bP[\cE_3]\leq \mu n^{2k-1}\cdot \left(\frac{2tmk}{n}\right)^{2k-1}\leq \frac{1}{m}$, as required.

\subsection*{Bounding $\bP[\cE_4]$.} Let $S\in \binom{V(G)}{j}$ be some $j$-set in $V(G)$ and define $E_G(S)-S:=\{e\setminus S:e\in E_G(S)\}\subset \binom{V(G)}{k-j}$ to be the collection of subsets that form an edge with $S$. Moreover, let $N'(S)\subseteq E_G(S)-S$ be the subcollection of sets $A$ such that $A$ is spread and no vertex in $A$ is close to a vertex in  $S\cup P_0$. Note that $N'(S)$ can be obtained from $E_G(S)-S$ by deleting at most $2tk^2n^{k-1-j}+n^{k-1-j}\cdot 2tk\leq \beta  n^{k-j}$ sets and so $|N'(S)|\geq (\delta_j^k(1)+\eps-\beta) n^{k-j}$ by our minimum $j$-degree condition on $G$ in the statement of the lemma. 

Let now $X(S):=|\{A\in N'(S):A\subset V(\cP_p)\}|$ be the random variable that counts the number of sets in $N'(S)$  whose vertices are contained in the random set $V(\cP_p)$.   The fact that $N'(S)$ does not contain any close pairs implies that vertices of sets in $N'(S)$ appear in $V(\cP_p)$ independently. Due to this as well as our lower bound on $N'(S)$ and the lower bound in \eqref{eq:single vx prob}, we have that 
 \begin{linenomath}
\begin{equation} \label{eq:lambda'}
\lambda':=\bE[X(S)]\geq (\delta_j^k(1)+\eps-\beta) n^{k-j}\cdot \left(\frac{(1-\beta)M}{n}\right)^{k-j}=(\delta_j^k(1)+7\eps/8)M^{k-j},\end{equation}
\end{linenomath}
using here that $\beta \ll \eps$. Let $\Gamma$ be the auxiliary graph with vertex set $N'(S)$ such that two sets $A,A'\in N'(S)$ are adjacent in $\Gamma$ whenever there is a pair of close vertices $u,u'$ with $u\in A$ and $u'\in A'$. Then $\Gamma$ is a strong dependency graph for the indicator random variables $\{A\subset V(\cP_p)\}_{A\in N'(S)}$ and we are in the setting of Suen's inequality (Lemma~\ref{lem:suen}). For $A,A'\in N(G')$ let $s(A,A')$ denote the size of a largest spread subset of $A\cup A'$ and note that when $AA'\in E(\Gamma)$ (equivalently, $A\sim A'$) one necessarily has $k-j\leq s(A,A')\leq 2(k-j)-1$. Then
in the notation of Lemma~\ref{lem:suen}, we have that 
\[\Delta_{X(S)}\leq \sum_{s=k-j}^{2(k-j)-1}\sum_{A\sim A':\\ s(A,A')=s}\bP[A\cup A'\subset V(\cP_p)]\leq \sum_{s=k-j}^{2(k-j)-1} n^s\cdot(4tk)^{2k}\left(\frac{2tmk}{n}\right)^{s}\leq (4tk)^{6k}m^{2(k-j)-1}, \]
using \eqref{eq:single vx prob}. Similarly we have that  
$\delta_{X(S)}\leq 2\cdot(4tk)\cdot n^{k-1-j}\left({2tmk}/{n}\right)^{k-j}\leq 1,$
using that sets in $N'(S)$ do not contain pairs of close vertices and the fact that $1/n\ll 1/m,1/t,1/k$ \eqref{eq:hierarchy}. Therefore appealing to  Lemma~\ref{lem:suen},  and using \eqref{eq:lambda'} we get that 
\begin{linenomath}
\begin{align}
\bP\left[X(S)< \left(\delta_j^k(1)+\frac{3\eps}{4}\right) M^{k-j}\right]&\leq\bP\left[X(S)\leq \lambda'-\frac{\eps}{8} M^{k-j}\right] \nonumber \\ &\leq \exp \left(-\frac{\eps^2 M^{2(k-j)}}{512(4tk)^{6k}m^{2(k-j)-1}}\right) \leq \exp \left(-\sqrt{m}\right),
\end{align}
\end{linenomath}
where in the last inequality we used  $1/m\ll 1/t,\eps,1/k$.

Finally then we need to sum over all possible $S\in \binom{V(G)}{j}$. For any such $S$, let $s(S)$ be the maximal size of a  subset of $S$ which is spread and contains no vertices close to $P_0$. Then $\bP[S\subset V(\cP_p)]\leq (2tmk/n)^{s(S)}$ using \eqref{eq:single vx prob} and the fact that vertices in such a maximal subset appear independently. Moreover, as sets in $N'(S)$ do not contain vertices close to $S$, we have that the events $\{S\subset V(\cP_p)\}$ and $\left\{X(S)< \left(\delta_j^k(1)+3\eps/4\right) M^{k-j}\right\}$ are independent. Thus  by a union bound
\begin{linenomath}
\begin{align*}
\bP[\cE_4] &\leq \sum_{S\in \binom{V(G)}{j}} \bP[S\subset V(\cP_p)]\cdot \bP\left[X(S)< \left(\delta_j^k(1)+3\eps/4\right) M^{k-j}\right]\\
&\leq \sum^{j}_{s=0} \sum_{\{S: s(S)=s\}} \left(\frac{2tmk}{n}\right)^s\exp(-\sqrt{m}) \\
&\leq \sum^{j}_{s=0} n^s\cdot ((s+tk)(4tk))^{j-s}\cdot \left(\frac{2tmk}{n}\right)^s\cdot \exp(-\sqrt{m}) \\ &\leq \frac{1}{m},
\end{align*}
\end{linenomath}
where the second summand in the penultimate line accounts for choosing vertices in $S$ that are close to a vertex in $P_0$ or a vertex in the maximal spread subset of $S$ with size $s(S)$. 

\subsection*{Bounding $\bP[\cE_5]$.}
The expected number of  vertices in $\cP_p$ that are close to some vertex in $V(P_0)$ is at most 
\[(3t+1)k\cdot \left(\frac{2tmk}{n}\right)\leq \frac{1}{2m},\]
by \eqref{eq:single vx prob}. Moreover, again using the upper bound of \eqref{eq:single vx prob}, we have that the expected number of pairs of vertices that are close and appear in distinct paths in $\cP_p$  is at most 
\[ n\cdot 4tk\cdot \left(\frac{2tmk}{n}\right)^2\leq \frac{1}{2m}. \]
Therefore by Markov's inequality, the probability that there is some pair of close vertices that appear in distinct paths in $\cP_*$ is at most $1/m$, as required.
\end{proof}

Next we show that any splitting as in Lemma~\ref{lem:random split} admits a viable partition. We will first consider a random partition and show that it has some of the desired properties, for example satisfying the degree condition~\ref{via:deg} of Definition~\ref{def:viable}. We will then  utilise the directed version of Dirac's theorem for Hamilton cycles in digraphs \cite{ghouilahouri1960condition} (see Theorem~\ref{thm:dirac_di} below)  to dictate some small adjustments to this random partition that will result in a viable partition.  

Let us first recall some relevant notation and terminology for digraphs. A \emph{digraph} $\vec{D}$ has a vertex set $V(\vec{D})$ and edge set $\vec{E}(\vec{D})$ consisting of ordered pairs $\vec{uv}$ of vertices $u\neq v\in V(\vec{D})$. The \emph{out-degree} of a vertex $u\in V(\vec{D})$ is $d^+(u):=|\{\vec{uv}\in \vec{E}(\vec{D}):v\in V(\vec{D})\}|$ and the \emph{in-degree} of $u\in V(\vec{D})$ is   $d^-(u):=|\{\vec{vu}\in \vec{E}(\vec{D}):v\in V(\vec{D})\}|$. The \emph{minimum out-degree} is $\delta^+(\vec{D}):=\min\{d^+(u):u\in V(\vec{D})\}$ and the \emph{minimum in-degree} $\delta^-(\vec{D})$ is defined analogously. 
A \emph{Hamilton dicycle} $\vec{H}$ in an $n'$-vertex digraph $\vec{D}$ is a subdigraph of $\vec{D}$ defined by the edge set $E(\vec{H})=\{\vec{u_iu}_{i+1}:i=0,\ldots, n'-1\}$ for some labelling of the vertices of $\vec{D}$ as $V(\vec{D})=\{u_0,\ldots, u_{n'-1}\}$ (and taking $u_{n'}:=u_0$).   
Ghouila-Houri \cite{ghouilahouri1960condition} proved the following natural generalisation of Dirac's theorem to digraphs.

\begin{thm}[Dirac's theorem for digraphs \cite{ghouilahouri1960condition}] \label{thm:dirac_di}
If $\vec{D}$ is an $n'$-vertex digraph with $\delta^-(\vec{D}),\delta^+(\vec{D})\geq n'/2$, then $\vec{D}$ contains a Hamilton dicycle. 
\end{thm}

We are now in a position to prove the second part of Lemma~\ref{lem:many suitable}, by finding a viable partition of our suitable splitting. 

\begin{lem} \label{lem:partition}
For a loose Hamilton cycle $H$ in $G$ and a path $P_0\subset H$ of length $t$, suppose that $\cP$ is  a $t$-balanced $m$-splitting of $H$ suitable for $P_0$ and $\delta_j(G[V(\cP)])\geq (\delta_j^k(1)+\tfrac{3\eps}{4})M^{k-j}$. Then $V(\cP)$ admits a viable partition.
\end{lem}
\begin{proof}
Let $\cP=\{P_0,\ldots,P_{m-1}\}$ and for each $i=0,\ldots,m-1$, label the endvertices of the path $P_i$ as $V^E(P_i):= \{u^i_a,u^i_b\}$ in such a way that if we start at some edge $e\in E(H)\setminus E(P_i)$ and traverse $H$ according to the implicit orientated cyclic order on $E(H)$, we meet the edge of $P_i$ containing $u^i_a$ before meeting the edge of $P_i$ containing $u^i_b$. Further, we fix $U_a:=\{u^i_a:i=0,\ldots,m-1\}$ and $U_b:=\{u^i_b:i=0,\ldots,m-1\}$ so that $V^E(\cP)=U_a\cup U_b$. Finally we follow the convention that $u^m_a:=u^0_a$, $u^m_b:=u^0_b$, $u^{-1}_a:=u^{m-1}_a$ and $u^{-1}_b:=u^{m-1}_b$.

We proceed with the following claim.

 \begin{clm} \label{clm:random part} There exists a transverse partition $\cY=\{Y_1,\ldots,Y_{T}\}$ of $V(\cP)$ with respect to $\cP$ (with parts of size $m$) such that for each $h\in [T]$ we have that 
\begin{enumerate}[label=(\Alph*)]
\item \label{vi:Ub} $|Y_h\cap U_b|=\tilde{m}$;
\item \label{vi:Ua} $|Y_h\cap U_a|\leq \beta m$; and 
\item \label{vi: better deg} $\deg_G(S;Y_h)\geq (\delta_j^k(1)+\tfrac{5\eps}{8})m^{k-j}$ for all $S\in \binom{V(\cP)}{j}$.
\end{enumerate} 
\end{clm}
\begin{claimproof}
The claim is proven by considering a uniformly random transverse partition. 
We will show that with probability at least $1/(4\sqrt{\tilde{m}})^{t(k-1)}$ property~\ref{vi:Ub} holds for all $Y_h$. On the other hand,  we will show that the probability that~\ref{vi:Ua} or~\ref{vi: better deg} fails is at most $e^{-\sqrt{m}}$. Hence  using that $\tilde{m}=t^{100k}$ and $m\geq \tilde{m}$ by \eqref{eq:m}, 
with positive probability the conditions are all satisfied and  such a transverse partition does indeed exist.

First let us lower bound the probability of property~\ref{vi:Ub} occurring for all $Y_h$. Note that we can generate the uniformly random partition $\cY$ via the following random process. We first choose $Y_1$ by taking a uniformly random vertex from each of the $V(P_i)$ with $i=0,\ldots, m-1$. Next, we choose $Y_2$ by taking a uniformly random vertex from each of  the remaining sets $V(P_i)\setminus Y_1$ with $i=0,\ldots, m-1$ and so on. At the last step our random choice for $Y_{t(k-1)}$ also defines $Y_{T}$ as the only remaining vertex from each of the $V(P_i)$ with $i=0,\ldots, m-1$. Now we have that
 \begin{linenomath}
\begin{equation} \label{eq:prod}
\bP[\ref{vi:Ub}]=\prod_{h=1}^{t(k-1)}\bP\left[|Y_h\cap U_b|=\tilde{m} \hspace{1mm} \big| \hspace{1mm}|Y_{h'}\cap U_b|=\tilde{m} \mbox{ for all }h'<h \right],
\end{equation}
\end{linenomath}
noting again that $Y_{T}$ is fully determined by the previous choices of $Y_{h'}$ and if each of the choices satisfy $|Y_{h'}\cap U_b|=\tilde{m}$ then certainly we will have $|Y_{T}\cap U_b|=\tilde{m}$. Now for a given $h\in[t(k-1)]$, given that  $|Y_{h'}\cap U_b|=\tilde{m}$ for all $h'<h$, we have that there exactly $m-(h-1)\tilde{m}=\tilde{m}\cdot(t(k-1)+2-h)$ indices $i\in\{0,\ldots,m-1\}$ such that $u_b^i$ is available for $Y_h$ and not used in previous $Y_{h'}$. Each of these vertices is chosen for $Y_h$ with probability exactly $1/(t(k-1)+2-h)$. Thus in expectation we have that $|Y_h\cap U_b|=\tilde{m}$ and as the the choices are independent for the different $i\in \{0,\ldots,m-1\}$ we have that $|Y_h\cap U_b|$ is binomially distributed. Using Lemma~\ref{lem:hit expected} to lower bound each of the terms in the product in \eqref{eq:prod} thus gives our required lower bound on the probability. 

We now turn to~\ref{vi:Ua}. Fixing a $h\in [T]$  we have that each vertex in $U_a$ lies in $Y_h$ with probability $1/T$. Hence $\bE[|U_a\cap Y_h|]=m/T=\tilde{m}$. Moreover the events $\{u\in Y_h\}$ are independent for different $u\in U_a$. Therefore we are in the setting of Chernoff's inequality (Theorem~\ref{thm:chernoff}) and as $1/t\ll \beta$, we get that $\bP[|U_a\cap Y_h|\geq \beta m/4]\leq \exp(-\beta m/4)$. Taking a union bound over the choices of $h\in [T]$, we get that with probability at least $1- \exp(-\sqrt{m})$ all the constraints on the sizes of these intersections are satisfied. 

 Finally, we address the degree condition~\ref{vi: better deg}. For $S\in \binom{V(\cP)}{j}$, define
\[\hat N(S):=\left\{A=e\setminus S: e\in E_{G[V(\cP)]}(S) \mbox{ and } A \mbox{ is transverse with respect to }\cP \right\}\]
As $\hat N(S)$ is obtained from the collection of all sets of the form $e\setminus S$ with $e\in E_{G[V(\cP)]}(S)$ by deleting at most $M^{k-j-1}\cdot(tk)\leq \beta M^{k-j}$ sets, we have that $|\hat N(S)|\geq (\delta_j^k(1)+\tfrac{3\eps}{4}-\beta)M^{k-j}$ for all $S$, using the minimum degree condition in the statement of Lemma~\ref{lem:partition}.  As there are at most $M^{j}$ $j$-sets in $V(\cP)$ and $T$ different $Y_h$, using a union bound, it suffices to show that for any fixed $S\in \binom{V(\cP)}{j}$ and any $h\in [T]$, we have that $\bP[\hat Z_h(S)\leq (\delta_j^k(1)+\tfrac{5\eps}{8})m^{k-j}]\leq \exp(-{m^{3/4}})$, where $\hat Z_h(S)=|\{A\in \hat N(S):A\subseteq Y_h\}|$ counts the number of sets in $\hat N(S)$ that fall into $Y_h$. Now for any such $S$ and $h\in[ T]$ we have that $A\subset Y_h$ with probability $(1/T)^{k-j}=(m/M)^{k-j}$, as it is transverse, and so using our lower bound on $|\hat N(S)|$, we get that $\hat \lambda:=\bE[\hat Z_h(S)]\geq (\delta_j^k(1)+\tfrac{3\eps}{4}-\beta)m^{k-j}$. An application of Suen's inequality (Lemma~\ref{lem:suen}) will then give exponential concentration to this mean and imply that indeed $\bP[\hat Z_h(S)\leq (\delta_j^k(1)+\tfrac{5\eps}{8})m^{k-j}]\leq \exp(-{m^{3/4}})$, as required. As the calculations for $\Delta_{\hat Z(S)}$ and $\delta_{\hat Z(S)}$ are very similar to those done in bounding $\bP[\cE_4]$ in the proof of Lemma~\ref{lem:random split}, we skip the details here. 
\end{claimproof}

\begin{figure}[ht] \begin{center}

    \includegraphics[width=14cm, height=6.8cm]{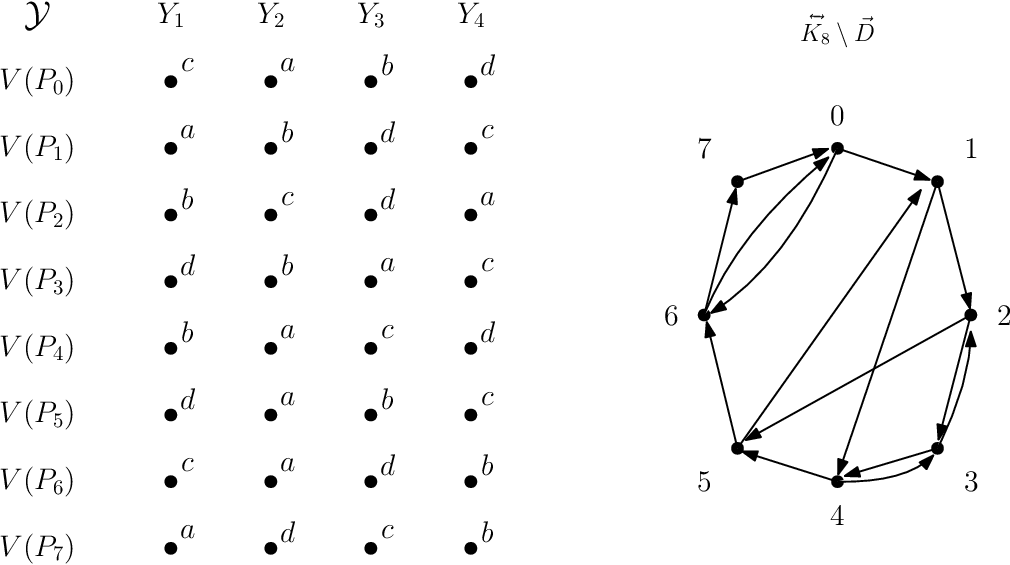}
  \caption{On the left, an example of a transverse partition $\cY=\{Y_1,Y_2,Y_3,Y_4\}$ with paths $P_0,\ldots,P_7$ (which are just single edges of a $4$-uniform hypergraph in this case) such that $V(P_i)=\{u^i_a, u^i_b, u^i_c,u^i_d\}$ (denoted in the figure simply as $a,b,c,d$ in each row) and $u^i_a$ and $u^i_b$ are the endpoints of $P_i$. 
  On the right, the digraph $\dvec{K_8}\setminus \vec{D}$ where $\dvec{K_8}$ is the complete digraph on $8$ vertices with all edges $\vec{ii'}$ with $0\leq i<i'\leq 8$, and $\vec{D}$ is the auxiliary digraph defined by $\cY$.}\label{fig:partition1}   
 \end{center}
      \end{figure}

\vspace{2mm}

Claim~\ref{clm:random part} puts us in good stead to prove the lemma. Indeed, the partition $\cY$ is transverse with respect to $\cP$ and condition~\ref{via:deg} of Definition~\ref{def:viable} is already satisfied (with some room to spare) due to condition~\ref{vi: better deg}. However, we are missing condition~\ref{via:reroute} of being viable and need to give a rerouting which is compatible with the partition. Claim~\ref{clm:random part} gives us some help to do this as we already have that each part $Y_h$ contains the correct number of vertices from $U_b$ that we need for a rerouting. Hence we will leave these $U_b$ vertices fixed and move the $U_a$ vertices, slightly modifying the partition $\cY$ in the process, in order to get some viable partition $\cX$.


Now to find the reshuffling, we appeal to an auxiliary digraph $\vec{D}$ with vertex set $V(\vec{D})=\{0,\ldots, m-1\}$ and for $i\neq i'$ we have $\vec{ii'}\in \vec{E}(\vec{D})$ if and only if $u^i_b$ and $u^{i'-1}_b$ lie in different sets $Y_h$. An example of (the complement of) this digraph  is given on the right side of Figure~\ref{fig:partition1}. Note that, by definition, $\vec{ii}_+\notin \vec E(\vec{D})$ where $i_+:=i+1 \mod m$, for all $i=0,\ldots,m-1$. 
Due to property~\ref{vi:Ub} given by Claim~\ref{clm:random part}, we have that $\delta^-(\vec{D}),\delta^+(\vec{D})\geq m-1-\tilde{m}\geq m/2$ and therefore by Theorem~\ref{thm:dirac_di}, $\vec{D}$ contains some Hamilton dicycle $\vec{H}$. See Figure~\ref{fig:Hamdi} for an example of such a Hamilton  dicycle  for the partition from Figure~\ref{fig:partition1}.

Letting $\cR:=\{\{u_a^i,u^{i'-1}_b\}: \vec{ii'}\in \vec{E}(\vec{H})\}$, we have that $\cR$ is a \emph{rerouting} (Definition~\ref{def:reroute}).  Indeed, any path in $H[V(G)\setminus V^I(\cP)]$ has endpoints $u_b^{i-1}$ and $u_a^i$ for some $i\in \{0,\ldots, m-1\}$ and identifying a pair of vertices in $\cR$ 
results in joining two of these paths by identifying an endpoint from each. The fact that $\vec{H}$ is a Hamilton dicycle guarantees that identifying all the pairs in $\cR$ results in a single loose cycle. See Figure~\ref{fig:partitions} for the pairs $\cR$ defined by the Hamilton dicycle in Figure~\ref{fig:Hamdi} for the partition $\cY$ in Figure~\ref{fig:partition1}.

\begin{figure}[ht]
 \begin{center}
    \includegraphics[width=5cm, height=5cm]{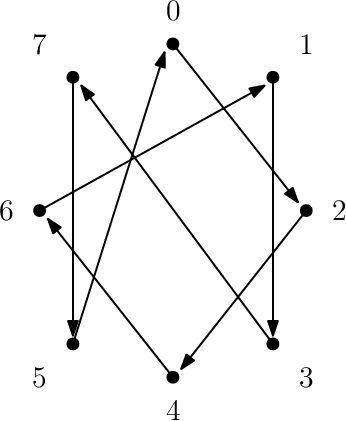}
  \caption{An example of Hamilton dicycle $\vec{H}$ in the digraph $\vec{D}$ given by the partition $\cY$ given in Figure~\ref{fig:partition1}.}\label{fig:Hamdi}
    \end{center}
      \end{figure}

Now for $i\in \{0,\ldots, m-1\}$, let $h_i\in [T]$ be the index such that $u_b^{i'-1}\in Y_{h_i}$ where $i'\in\{0,\ldots,m-1\}$ is such that $\vec{ii'}\in \vec{E}(\vec{H})$. Further, define $w_i\in V(\cP)$ to be the unique vertex lying in the intersection of $V(P_i)$ and $Y_{h_i}$. Note that due to how we defined the auxiliary digraph $\vec{D}$, we have that $w_i\neq u_i^b$ as we guaranteed that for each edge $\vec{ii'}\in \vec{E}(\vec{H})\subset \vec{E}(\vec{D})$ we have that $u^i_b$ and  $u_b^{i'-1}$ lie in different parts of $\cY$. It may be the case however that $u_a^i\in Y_{h_i}$ and thus $w_i=u_a^i$. Examples are given in Figure~\ref{fig:partitions}.

We now define $\cX:=\{X_1,\ldots, X_{T}\}$ to be the partition of $V(\cP)$  obtained from $\cY$ by swapping all the pairs $\{\{u_a^i,w_i\}:i=0,\ldots,m-1\}$. 
Formally, we have that for $h\in [T]$, 
\[X_h:=\left(Y_h\setminus (U_a\cup\{w_{h_i}:i\in I_h\})\right)\cup \{w_i:i\in \{0,\ldots,m-1\} \mbox{ with }u_a^i\in Y_h\}\cup\{u_a^i:i\in I_h\},\]
where $I_h:=\{i\in \{0,\ldots,m-1\}:h_i=h\}$. See Figure~\ref{fig:partitions} for an example.

\begin{figure}[ht]
 \begin{center}   
    \includegraphics[width=16cm, height=6.8cm]{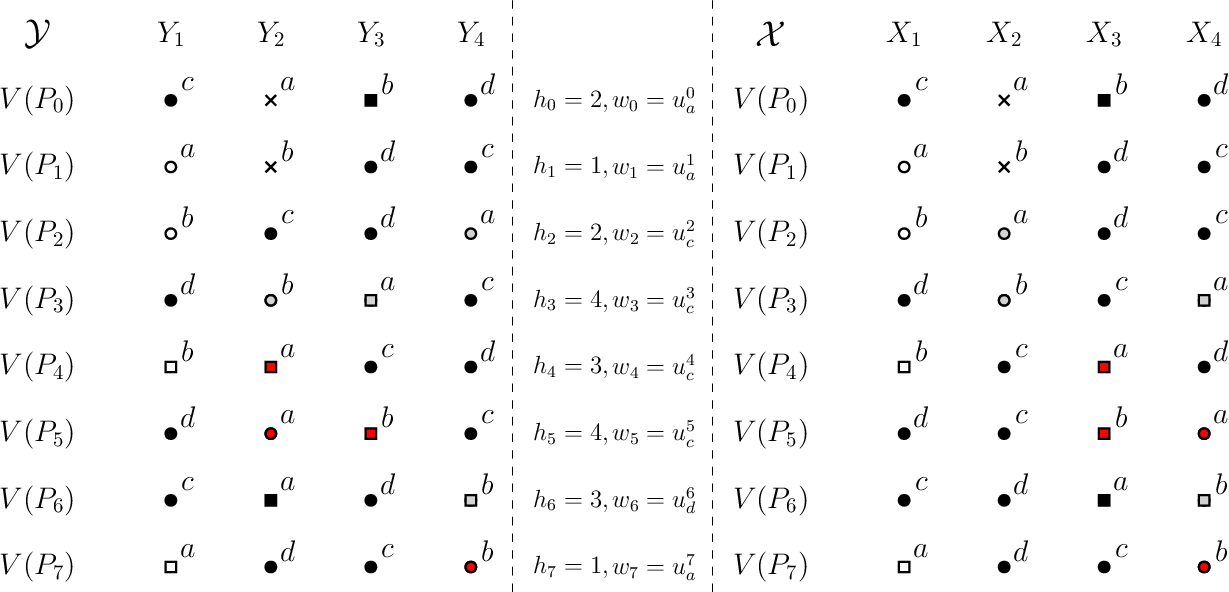}
  \caption{On the left, the pairs $\cR$ are identified for the partition $\cY$ according to the Hamilton dicycle $\vec{H}$ from Figure~\ref{fig:Hamdi}. The values of $h_i$ and the vertices $w_i$ are then given and the resulting partition $\cX$. 
  }   \label{fig:partitions}
 \end{center}
      \end{figure}

Note that whenever we remove a vertex in $Y_h\cap V(P_i)$ for some $i\in\{0,\ldots,m-1\}$ it is replaced by some other vertex in $Y_h\cap V(P_i)$ and so $\cX$ is again a transverse partition of $V(\cP)$ into parts of size $m$. Moreover, due to condition~\ref{vi:Ub} of Claim~\ref{clm:random part}, the fact that $w_i\notin U_b$ for all $i\in\{0,\ldots,m-1\}$ and how we defined  the partition $\cX$, we have that each $X_h$ contains exactly $\tilde{m}$ pairs of the rerouting $\cR$.  Thus condition~\ref{via:reroute} of being viable (Definition~\ref{def:viable}) is satisfied. Note also that for each $h\in [T]$, we have that
\[|Y_h\setminus X_h|\leq |Y_h\cap U_a|+ |I_h|\leq \beta m+\tilde{m}\leq \eps m/8, \]
using properties~\ref{vi:Ub} and~\ref{vi:Ua} of Claim~\ref{clm:random part} and the fact that any index in $I_h$ corresponds to a vertex in $U_b\cap Y_h$. Thus for any $S\subset X_h$ with $|S|=j$, there are at most $\eps m^{k-j}/8$ sets in $\binom{Y_h}{k-j}\setminus \binom{X_h}{k-j}$. Hence the $j$-degree of $S$ relative to $Y_h$ is close to its $j$-degree in $X_h$. Thus, due to property~\ref{vi: better deg} of Claim~\ref{clm:random part}, we get that for each $h\in[T]$, we have $\delta_j(G[X_h])\geq (\delta_j^k(1)+\tfrac{\eps}{2})m^{k-j}$, verifying  property~\ref{via:deg} of being viable (Definition~\ref{def:viable}) and completing the proof of the lemma. 
\end{proof}

With Lemmas~\ref{lem:random split} and~\ref{lem:partition}, Lemma~\ref{lem:many suitable} follows easily. 

\begin{proof}[Proof of Lemma~\ref{lem:many suitable}]
Recall that $H$ is a loose Hamilton cycle and $P_0\subset H$. Let $\bc{P}$ denote the collection of $t$-balanced $m$-splittings $\cP$ of $H$ that are suitable for $P_0\in \cP$ and such that $V(\cP)$ admits a viable partition.  Moreover, let  $\cP_*$ be the output by the random procedure in Lemma~\ref{lem:random split}. Then by Lemmas~\ref{lem:random split} and~\ref{lem:partition}, we have that 
\[1/m\leq \bP[\cP_*\in \bc{P}]=\sum_{i=1}^{n+1}\bP\left[\cP_*\in \bc{P}| \hspace{1mm} |\cP_*|=i\right] \bP[|\cP_*|=i] \leq \bP\left[\cP_*\in \bc{P}| \hspace{1mm} |\cP_*|=m\right], \]
using that $\bP\left[\cP_*\in \bc{P}| \hspace{1mm} |\cP_*|=i\right]=0$ for all $i\neq m$ in the last inequality. Given that $|\cP_*|=m$, each of the possible $\binom{n/(k-1)}{m-1}$ outcomes of $\cP_*$ (given by a choice of $E_p\subset E(H)$ with $|E_p|=m-1$), are equally likely. Hence 
\[|\bc{P}|= \bP\left[\cP_*\in \bc{P}| \hspace{1mm} |\cP_*|=m\right]\binom{\tfrac{n}{k-1}}{m-1}\geq \frac{1}{m}\binom{\tfrac{n}{k-1}}{m-1}\geq \gamma n^{m-1},\]
as required, using here that $\gamma \ll 1/m, 1/k$. 
\end{proof}

\section{Concluding Remarks} \label{sec:conclude}
In this paper, we made progress towards the meta-conjecture of Coulson, Keevash, Perarnau and Yepremyan \cite{ckpy2020rainbow} that whenever one is above the extremal threshold for a given spanning structure, one can find a rainbow copy of that spanning structure in any suitably bounded colouring. Indeed, Theorem \ref{thm:main} confirms this phenomenon for loose Hamilton cycles in hypergraphs. We also elucidate what constitutes a \emph{suitable} bound on the colouring,  with the global bound \ref{item:Global} in Theorem \ref{thm:main} capturing the optimal qualitative conditions in the setting of loose Hamilton cycles. As mentioned in the introduction, for $F$-factors in hypergraphs, the theorem of Coulson, Keevash, Perarnau and Yepremyan \cite{ckpy2020rainbow} requires an extra \emph{local} bound that every $(k-1)$-set is in $o(n)$ edges of the same colour. This is necessary in some cases. Indeed, with $k\geq 3$ consider $G$ to be the complete $k$-graph and assign a unique colour to each $(k-1)$-set of vertices of $G$. If we then colour the edges of $G$ by ordering its vertices and colouring each edge by the colour corresponding to the first $(k-1)$-set of vertices in the edge, we have that this colouring is globally but not locally bounded. Moreover, if $F$ is such that every $(k-1)$-set of vertices of $F$ is contained in at least two edges of $F$ (for example, if $F$ is a clique on at least $k+1$ vertices), then one can see that any $F$-factor in $G$ is \emph{not} rainbow. This justifies the necessity of the local bound in general. However, as mentioned in the introduction, one can use our proof, in particular our method in Lemma \ref{lem:path tilings} using a partition to avoid edges of the bad graph $B$, to remove the need for a local bound when $F$ is \emph{linear}, that is, when no pair of edges of $F$ intersect in more than one vertex. It remains to fully determine optimal ``suitable'' conditions on bounds for \emph{all} $F$-factors.

 Even more pertinent is the meta-conjecture for other types of Hamilton cycles in hypergraphs, where we lack any sort of positive result. 
Focusing in on the other extreme to loose Hamilton cycles in hypergraphs, 
 for tight Hamilton cycles the only result known is due to Dudek, Frieze and Ruci\'nski \cite{dudek2012rainbow} who showed the existence of rainbow tight cycles in the complete $k$-uniform hypergraph and under a strong global bounded condition of having $o(n)$ edges of each colour. It would be interesting to generalise this  result to Dirac hypergraphs and establish an analogue of Theorem \ref{thm:main} in the context of tight cycles (or indeed any Hamilton $\ell$-cycle with $\ell\geq 2$). We remark that the proof scheme used here which appeals to 
 switchings and double counting, seems to break down for these cases as tight Hamilton cycles are, in a sense, \emph{more} connected and so there is less flexibility to locally alter a Hamilton cycle to obtain (enough) new ones. Ultimately, it would be interesting to establish optimal conditions on the boundedness needed to guarantee the existence of tight Hamilton cycles (and the other $\ell$-cycles) in Dirac graphs, as we have for loose cycles. This was already raised explicitly in \cite{antoniuk2023properly}. Even in the complete case, this remains mysterious and Dudek, Frieze an Ruci\'nski \cite{dudek2012rainbow} asked if  some strong global bound (or a local bound as above) is in fact necessary. That is, they asked if there is some colouring of the complete $k$-graph with $\Theta(n)$ edges of each colour and no rainbow tight Hamilton cycle. They also suggested that they believe the answer should be positive. Alternatively, it could be that a weak global bound of $o(n^{k-1})$ edges of each colour (as in Theorem \ref{thm:main}) already suffices to force a rainbow tight Hamilton cycle. This remains open but the following example shows that, at least in some range of minimum degree, bounding the colouring further than the weak global bound  \emph{is} necessary. 

\begin{figure}[ht] 
 \begin{center}                
    \includegraphics[width=6cm, height=6cm]{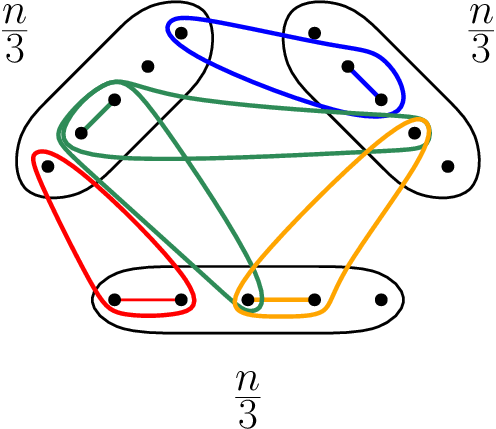}
\caption{The construction for Example \ref{example:tight} in the case that $n\in 3 \NN$.} \label{fig:exampleyy}
      \end{center}
\end{figure}

 \begin{exmp} \label{example:tight}
     For any $0<\alpha<2/3$ and $n$ sufficiently large, there exists an $n$-vertex $3$-graph $G$  with $\delta_2(G)\geq \alpha n$ and a colouring  of $G$ with no colour appearing more than $n$ times such that there is no rainbow tight Hamilton cycle.
 \end{exmp}
 \begin{proof}
     We first define the graph $G$ and the colouring, see Figure \ref{fig:exampleyy}. Partition the vertex set $V$ of $G$ into three vertex subsets $V_1,V_2,V_3$ such that $|V_i|\in \{\floor{n/3}\ceil{n/3}\}$ for each $i\in [3]$. For every pair of vertices $\{u,v\}\in \binom{V}{2}$ such that $\{u,v\}\subset V_i$ for some $i\in [3]$, define a unique colour $c(u,v)$. Now the edges of $G$ are precisely the sets $S=\{u,v,w\}$ such that $|S\cap V_{i'}|=2$ for some $i'\in [3]$ and we define the colouring of $S\in E(G)$ to be $c(S\cap V_{i'})$. 
     The fact that $\delta_2(G)\geq \alpha n$ follows from the fact that $\alpha<2/3$ and $n$ is sufficiently large. It is also clear that each colour appears at most $2\ceil{n/3}\leq n$ times. To see that $G$ does not contain a rainbow tight Hamilton cycle, suppose for a contradiction that it does and let $v_1,v_2\ldots,v_n$ be a cyclic ordering of $V$ corresponding to the cycle. Moreover, without loss of generality, suppose that $v_1\in V_1$ and $v_2\in V_2$. Let $t\in \NN$ be the least integer such that $v_t\in V_3$ and consider $v_{t-2},v_{t-1}\in V_1\cup V_2$. We must have that there is some $i_*\in [2]$ such that $v_{t-2},v_{t-1}\in V_{i_*}$. Indeed otherwise, as the ordering corresponds to a tight Hamilton cycle, there would be an edge $\{v_{t-2},v_{t-1},v_t\}\in E(G)$ with one vertex in each part, which is not possible. Now  $v_{t-3}\notin V_{i_*}$ as no edge is completely contained in $V_{i_*}$ and so we must have that $v_{t-3}\in V_{3-i_*}$.  Therefore both $\{v_{t-3},v_{t-2},v_{t-1}\}$ and $\{v_{t-2},v_{t-1},v_t\}$ lie on the Hamilton cycle and have colour $c(\{v_{t-2},v_{t-1}\})$, meaning that the Hamilton cycle is not rainbow, a contradiction.
 \end{proof}

It was proven by R\"odl, Ruci\'nski and Szemer\'edi  \cite{RRSz06} that the codegree threshold for tight Hamitlon cycles  $\delta_2^3(2)$ (see Definition \ref{def:threshold}) is located at $1/2$. Therefore Example \ref{example:tight} shows that in the minimum degree range $1/2<\delta(G)<2/3$, there are $3$-graphs above the extremal threshold for  which an $o(n^2)$ global bound on the  colouring does not  force a rainbow Hamilton cycle and some extra local bound is necessary. In general, we believe that the optimal bounds on colourings for tight Hamilton cycles (and other Hamilton $\ell$-cycles also) may depend on the minimum degree of the Dirac hypergraph in question and it would be interesting to explore this further. 

Finally we remark that we have only focused on qualitative features of bounding colourings and one may also ask for the best possible quantitative bounds, for example trying to optimise the constant $\mu>0$ in the global bound \ref{item:Global} in Theorem \ref{thm:main}. This seems like a very challenging problem in general, even for the simplest case of perfect matchings in globally bounded colourings of complete bipartite graphs. Here, the constant of Erd\H{o}s and Spencer \cite{erdos1991lopsided} has been improved slightly by Bissacot, Fern{\'a}ndez, Procacci and Scoppola \cite{bissacot2011improvement} who showed that if no colour appears more than $27n/256$ times, then there is a rainbow perfect matching. From above, it is clear that a global bound of $n$ is necessary in order to have enough colours. It follows from a recent  construction of Pokrovskiy and Sudakov \cite{pokrovskiy2019counterexample} (see also the very recent improvement \cite{chakraborti2024almost}) that being globally $n$-bounded is not enough to guarantee a rainbow perfect matching, but it is still possible that a bound of the form $n-o(n)$ suffices. 

\bibliographystyle{abbrv}
\bibliography{biblio.bib}

\end{document}